\newcommand{\R}{\mathds R}

\newcommand{\Levi}{\mathfrak L}
\newcommand{\Ker}{\mathrm{Ker}}
\newcommand{\dd}{\mathrm d}
\newcommand{\Ddt}[1]{\tfrac{\mathrm D^{\hbox to 0pt{$\scriptscriptstyle#1$\hskip 0pt minus 1fil}}}{\dd t}}
\newcommand{\Id}{\mathrm{Id}}
\newcommand{\bil}[1]{{#1}^\mathrm{bil}}
\newcommand{\Lin}{\mathrm{Lin}}
\newcommand{\Bilin}{\mathrm{Bilin}}

\newcommand{\ad}{\mathrm{ad}}
\newcommand{\hor}{\mathrm{hor}}

\newcommand{\gl}{\mathrm{gl}}
\newcommand{\GL}{\mathrm{GL}}

\newcommand{\Dom}{\mathrm{Dom}}
\newcommand{\Gr}{\mathrm{Gr}}
\newcommand{\B}{\mathrm B}
\newcommand{\Vect}[1]{\underline{\mathfrak{Vec}}^{#1}}

\renewcommand{\contentsline}[3]{\csname new#1\endcsname{#2}{#3}}
\newcommand{\newchapter}[2]{\bigskip\hbox to \hsize{\vbox{\advance\hsize by -.5cm\baselineskip=12pt\parfillskip=0pt\leftskip=2cm\noindent\hskip -2cm #1\leaders\hbox{.}\hfil\hfil\par}$\,$#2\hfil}}
\newcommand{\newsection}[2]{\medskip\hbox to \hsize{\vbox{\advance\hsize by -.5cm\baselineskip=12pt\parfillskip=0pt\leftskip=2.5cm\noindent\hskip -2cm #1\leaders\hbox{.}\hfil\hfil\par}$\,$#2\hfil}}
\newcommand{\newsubsection}[2]{\medskip\hbox to \hsize{\vbox{\advance\hsize by -.5cm\baselineskip=12pt\parfillskip=0pt\leftskip=3.5cm\noindent\hskip -2cm #1\leaders\hbox{.}\hfil\hfil\par}$\,$#2\hfil}}

\documentclass[oneside,draft,11pt]{amsart}

\usepackage{times}
\usepackage{dsfont}
\usepackage[all]{xy}

\title[The single-leaf Frobenius Theorem]{The single-leaf Frobenius Theorem with Applications}
\author[P.\ Piccione]{Paolo Piccione}
\author[D.\ Tausk]{Daniel V.\ Tausk}

\address{Departamento de Matem\'atica,\hfill\break\indent  Universidade de S\~ao Paulo,
Brazil}
\email{piccione@ime.usp.br, tausk@ime.usp.br}
\urladdr{http://www.ime.usp.br/\~{}piccione,\hfill\break\phantom{URL: } http://www.ime.usp.br/\~{}tausk}

\subjclass[2000]{53B05, 53C05, 53C42, 55R25}

\keywords{Smooth distributions, Levi form, Frobenius theorem, affine connections, Levi--Civita connections,
Cartan--Ambrose--Hicks theorem}

\dedicatory{Dedicated to Prof.\ Serge Lang}

\date{October 25th, 2005}

\begin{document}


\theoremstyle{plain}\newtheorem{teo}{Theorem}[section]
\theoremstyle{plain}\newtheorem{prop}[teo]{Proposition}
\theoremstyle{plain}\newtheorem{lem}[teo]{Lemma}
\theoremstyle{plain}\newtheorem{cor}[teo]{Corollary}
\theoremstyle{definition}\newtheorem{defin}[teo]{Definition}
\theoremstyle{remark}\newtheorem{rem}[teo]{Remark}
\theoremstyle{plain} \newtheorem{assum}[teo]{Assumption}
\swapnumbers
\theoremstyle{definition}\newtheorem{example}{Example}[section]


\begin{abstract}
Using the notion of Levi form of a smooth distribution, we discuss
the local and the global problem of existence of \emph{one} horizontal
section of a smooth vector bundle endowed with a horizontal
distribution. The analysis will lead to the formulation of
a ``one-leaf'' analogue of the classical Frobenius integrability
theorem in elementary differential geometry. Several applications
of the result will be discussed. First, we will give a characterization
of symmetric connections arising as Levi--Civita connections
of semi-Riemannian metric tensors. Second, we will
prove a general version of the classical Cartan--Ambrose--Hicks
Theorem giving conditions on the existence of an affine map
with prescribed differential at one point between manifolds
endowed with connections.
\end{abstract}

\maketitle
\tableofcontents
\begin{section}{Introduction}
The central theme of the paper is the study of conditions for the existence
of \emph{one} integral leaf of (non integrable) smooth distributions satisfying
a given initial condition. The integrability condition given by Frobenius
theorem,  a very classical result in elementary Differential Geometry,
guarantees the existence of integral leaves with \emph{any} initial
condition. If on one hand such condition is very strong, on the other hand
the involutivity assumption in Frobenius theorem is very restrictive.
For instance, the integrability of the horizontal distribution of
a connection in a vector bundle is equivalent to the flatness of
the connection.

A measure of non integrability for a smooth distribution $\mathcal D$
on a manifold $E$ is provided
by the so-called {\em Levi form} $\mathfrak L^{\mathcal D}$ of $\mathcal D$; this is a skew-symmetric
bilinear tensor defined on the distribution, taking values in the quotient
$TE/\mathcal D$. For $x\in E$ and $v,w\in \mathcal D_x$, the value $\mathfrak L^{\mathcal D}_x(v,w)$
is given by the projection on $T_xE/\mathcal D_x$ of the Lie bracket $[X,Y]_x$,
where $X$ and $Y$ are arbitrary extensions of $v$ and $w$ respectively to $\mathcal D$-horizontal
vector fields.
If $\Sigma\subset E$ is an integral submanifold of $\mathcal D$, then the Levi form of $\mathcal D$
vanishes on the points of $\Sigma$. The first central observation that is
made in this paper is that, conversely, given an immersed submanifold $\Sigma$ of $E$ with
$T_{x_0}\Sigma=\mathcal D_{x_0}$ for some $x_0\in E$,
if $\Sigma$ is \emph{ruled} (in an appropriate sense) by curves tangent to $\mathcal D$, and if
$\mathfrak L^{\mathcal D}$ vanishes along $\Sigma$, then
$\Sigma$ is an integral submanifold of $\mathcal D$.
In particular, assume that $\mathcal D\subset TE$ is a horizontal
distribution of a vector bundle $\pi:E\to M$ over a manifold
$M$, and that $\Sigma$ is a local section of $\pi$ which is obtained
by parallel lifting of a family of curves on $M$ issuing from some
fixed point $x_0$. If the Levi form of $\mathcal D$ vanishes along
$\Sigma$, then $\Sigma$ is a parallel section of $\pi$ (Theorem~\ref{thm:Frobenius});
we call this result the \emph{(local) single leaf Frobenius theorem}.
In the real analytic case, a higher order version of this result
is given in Theorem~\ref{thm:higherorderFrobenius}; roughly speaking, the
higher order derivatives of the Levi form $\mathfrak L^{\mathcal D}$
are obtained from iterated Lie brackets of $\mathcal D$-horizontal vector fields.
The higher order single-leaf Frobenius theorem states that, in the real-analytic case, if
at some point $x_0$ of the manifold $E$ all the iterated brackets of
vector fields in $\mathcal D$ belong to $\mathcal D_{x_0}$, then
there exists an integral submanifold of $\mathcal D$ through $x_0$.

A global version of the single-leaf Frobenius theorem is discussed
in Theorem~\ref{thm:globalFrobenius}; here, the base manifold $M$ has to
be assumed simply-connected. Assume that a spray is given on $M$, for instance,
the geodesic spray of some Riemannian metric.
The existence of a global parallel
section of $\pi$ through a point $e_0$ with $\pi(e_0)=x_0\in M$
is guaranteed by the following condition: every piecewise solution
$\gamma:[a,b]\to M$ of $\mathcal S$ with $\gamma(a)=x_0$
should admit a parallel lifting $\widetilde\gamma:[a,b]\to E$
such that $\widetilde\gamma(a)=e_0$ and such that the Levi form
of $\mathcal D$ vanishes at the point
$\widetilde\gamma(b)$.

We also observe (Proposition~\ref{thm:Frobeniusrealanalytic}) that
in the real analytic case, every local parallel section defined
on a non empty open subset of a simply connected manifold $M$ extends
to a global parallel section.
\smallskip

A huge number of problems in Analysis and in Geometry can be cast into
the language of distributions and integral submanifolds. As an application
of the theory developed in this paper, we will consider two problems.
First, we will characterize those symmetric connections that are
Levi--Civita connections of some semi-Riemannian metric (alternatively, this
problem can be studied using holonomy theory, see \cite{KN}).
Second, we will prove a very general version of another classical result
in Differential Geometry, which is the Cartan--Ambrose--Hicks theorem (see \cite{Wolf}).
We will prove a necessary and sufficient condition for the existence
of an affine map between manifolds endowed with arbitrary connections.
\smallskip

Let us describe briefly these two results.

Consider the case of a distribution given by the horizontal
space of a connection $\nabla$ of a vector bundle $\pi:E\to M$.
For $\xi\in E$, set $m=\pi(\xi)\in M$ and $E_m=\pi^{-1}(m)$;
one can identify $\mathcal D_\xi$ with $T_mM$, and the quotient $T_\xi E/\mathcal D_\xi$ with the
vertical subspace $T_\xi(E_m)\cong E_m$. Then, the Levi form
$\mathfrak L^{\mathcal D}_\xi: T_mM\times T_mM\to E_m$ is given by
the curvature tensor of $\nabla$, up to a sign (Lemma~\ref{thm:Levihorizontal}).
In this case, the single leaf Frobenius theorem tells us that
a local parallel section of $\pi$ through some point $e_0\in E$ exists provided that along each
parallel lifting of a family of curves issuing from
$\pi(e_0)\in M$ the curvature tensor vanishes (Corollary~\ref{thm:corparalelo}). In the real analytic case,
the existence of a local parallel section through a point
$\xi\in E$ is equivalent to the vanishing of all the covariant
derivatives $\nabla^kR$, $k\ge0$, of the curvature tensor $R$ at
the point $m=\pi(\xi)$ (Proposition~\ref{thm:prophighorderhorizontal}).

Assume that the vector bundle $\pi:E\to M$ is endowed with a connection
$\nabla$, and denote  by $\bil\nabla$ the induced connection on $E^*\otimes E^*$.
If $g$ is a (local) section of $E^*\otimes E^*$, then vanishing of the
curvature tensor $\bil R$ of $\bil\nabla$ means that the bilinear map
$g\big(R(v,w)\cdot,\cdot)$ is anti-symmetric for all $v,w$ (formula \eqref{eq:Rbil}).
From this observation, we get the following result on the existence of
parallel metric tensor relatively to a given connection $\nabla$
on a manifold $M$:
given a nondegenerate (symmetric) bilinear form $g_0$ on $T_{m_0}M$, assume
that the tensor $g$ obtained from $g_0$ by $\nabla$-parallel transport
along a family of curves issuing from $m_0$ is such that $R$ is
$g$-anti-symmetric. Then, $g$ is $\nabla$-parallel (Proposition~\ref{thm:vemdemetrica}).
Similarly, in the real analytic case, if $\nabla^kR$ at $m_0$ is $g_0$-anti-symmetric for
all $k\ge 0$, then $g_0$ extends to a semi-Riemannian metric tensor whose
Levi--Civita connection is $\nabla$. These results have been used in \cite{Cordoba}
to obtain characterizations of left-invariant semi-Riemannian  Levi Civita connections
in Lie groups.

As another application of our theory, in Section~\ref{sec:affine} we will study
the problem of existence of an affine (i.e., connection preserving) map $f$ between two affine manifolds
$(M,\nabla^M)$ and $(N,\nabla^N)$, whose value $y_0\in N$ at some point $x_0\in M$ is given and whose differential
$\mathrm df(x_0):T_{x_0}M\to T_{y_0}N$ is prescribed. We prove a general
version of the classical Cartan--Ambrose--Hicks theorem (Theorem~\ref{thm:localCAH}
for the local result, Theorem~\ref{thm:CAH} for the global version),
giving a necessary and sufficient condition for the existence of such a map; here, the connections
$\nabla^M$ and $\nabla^N$ are not assumed to be symmetric, and no assumption is made on the dimension
of the manifolds $M$ and $N$, as well as on the linear map $\mathrm df(x_0)$.
The key observation here (Lemma~\ref{thm:affineparallel}) is that, considering the vector bundle
$E=\Lin(TM,TN)$ over the product $M\times N$, endowed with a natural connection
induced by $\nabla^M$ and $\nabla^N$ (see formula \eqref{eq:indconnLinTmTn}), then a smooth
map  $f:U\subset M\to N$ is an affine map if and only if the differential $\mathrm df$ is a
local parallel section of $E$ along the map $U\ni x\mapsto\big(x,f(x)\big)\in M\times N$.
When $M$ and $N$ are endowed with semi-Riemannian metrics and $\nabla^M$ and $\nabla^N$ are the
respective Levi--Civita connections, then our result gives a necessary and sufficient condition
for the existence of a totally geodesic immersion of $M$ in $N$.

The proof of the Cartan--Ambrose--Hicks theorem is obtained as an
application of the single-leaf Frobenius theorem, once the Levi form
of the horizontal distribution of the induced connection on $E$ is computed
(Lemma~\ref{thm:LemaLeviD}). The higher order version of this result (Theorem~\ref{thm:highCAH}) is particularly
interesting: in the real analytic case, a (local) affine map $f:U\subset M\to N$
with $f(x_0)=y_0$ and $\mathrm df(x_0)=\sigma$ exists if and only if $\sigma$
relates covariant derivatives of all order of curvature and torsion of
$\nabla^M$ and $\nabla^N$ at the points $x_0$ and $y_0$ respectively.

As a nice corollary of the higher order Cartan--Ambrose--Hicks theorem,
we get the following curious result (Corollary~\ref{thm:curiouscor}):
if $M$ is a real-analytic manifold endowed with a real-analytic connection $\nabla$, and let $x_0\in M$ be fixed;
there exists an affine symmetry around $x_0$ if and only if $\nabla^{(2r)}T_{x_0}=0$
and $\nabla^{(2r+1)}R_{x_0}=0$ for all $r\ge0$.

\smallskip

A certain effort has been made in order to make the presentation
of the material self-contained. For this reason, we have found
convenient to discuss, together with the original material, some
auxiliary topics needed for a more complete presentation of our
results. For instance, in Subsection~\ref{sub:sprays}, we
discuss and give the basic properties of the exponential
map of a spray (this is needed in our statement of the global
one-leaf Frobenius theorem). Similarly, in Appendix~\ref{app:crashcourse}
we develop the basic theory needed for making computations with covariant
derivatives, curvatures and torsions of connections on vector bundles
obtained by functorial constructions; this kind of computations is
heavily used throughout the paper.
Finally, in Appendix~\ref{sec:globalization} we discuss a globalization
principle in a very general setting of pre-sheafs on topological spaces.
Such principle is used in the proof of the global versions of the single-leaf
Frobenius theorem (see for instance the proofs of Theorem~\ref{thm:globalFrobenius}
and Proposition~\ref{thm:Frobeniusrealanalytic}). Typically, the globalization
principle is employed in the following manner: given a vector bundle
$\pi:E\to M$, a pre-sheaf $\mathfrak P$ is defined on $M$ by defining,
for all open subset $U\subset M$, $\mathfrak P(U)$ to be the set of
all sections $s:U\to E$ of $\pi$ satisfying some property (for instance, parallel
sections). For $V\subset U$, and $s\in\mathfrak P(U)$, the map
$\mathfrak P_{U,V}:\mathfrak P(U)\to\mathfrak P(V)$ is given by setting
$\mathfrak P_{U,V}(s)=s\vert_V$. In this context, the existence
of a global section of $\pi$ with the required property is equivalent
to the fact that the set $\mathfrak P(M)$ should be non empty.
The central result of Appendix~\ref{sec:globalization} (Proposition~\ref{thm:niceglobalization})
gives a sufficient condition for this, in terms of three properties
of pre-sheaves, called \emph{localization}, \emph{uniqueness} and
\emph{extension}.

\smallskip

\noindent\textbf{Dedicatory.} The proof of the single-leaf Frobenius theorem discussed
here has taken inspiration from the proof of the classical Frobenius theorem presented
in Serge Lang's  world famous book \cite{Lang}. Since the very beginning
of their mathematical careers, both authors have benefited very much
from this and from other beautiful books published by Prof.\ Lang.
We want to thank him by dedicating this paper to his memory.
\end{section}

\begin{section}{The Levi form and the ``single leaf Frobenius Theorem''}
\label{sec:leviformsingleleaf}
Recall that a {\em smooth distribution\/} $\mathcal D$ on a smooth manifold $E$ is a smooth vector subbundle
of the tangent bundle $TE$. For $x\in E$ we set $\mathcal D_x=T_xE\cap\mathcal D$, i.e.,
$\mathcal D_x$ is the fiber of the vector bundle $\mathcal D$ over $x$. A vector field $X$ on $E$
is called {\em horizontal\/} with respect to a distribution $\mathcal D$ (or simply {\em $\mathcal D$-horizontal})
if $X$ takes values in $\mathcal D$, i.e., if $X(x)\in\mathcal D_x$ for all $x\in E$. An immersed
submanifold $S$ of $E$ is called an {\em integral submanifold\/} for $\mathcal D$ if $T_xS=\mathcal D_x$,
for all $x\in S$. The distribution $\mathcal D$ is called {\em integrable\/} if through every point of $E$ passes
an integral submanifold for $\mathcal D$.

\subsection{The Levi form of a smooth distribution}
\label{sub:levi}
\begin{defin}\label{thm:defLevi}
Let $E$ be a smooth manifold and let $\mathcal D$ be a distribution on $E$. The {\em Levi form\/}
of $\mathcal D$ at a point $x\in E$ is the bilinear map:
\[\Levi^{\mathcal D}_x:\mathcal D_x\times\mathcal D_x\longrightarrow T_xE/\mathcal D_x\]
defined by $\Levi^{\mathcal D}_x(v,w)=[X,Y](x)+\mathcal D_x\in T_xE/\mathcal D_x$, where $X$ and $Y$ are $\mathcal D$-horizontal
smooth vector fields defined in an open neighborhood of $x$ in $E$ with $X(x)=v$ and $Y(x)=w$.
By $[X,Y]$ we denote the Lie bracket of the vector fields $X$ and $Y$.
\end{defin}

Below we show that the Levi form is well-defined, i.e., $[X,Y](x)+\mathcal D_x$ does not
depend on the choice of the $\mathcal D$-horizontal vector fields $X$ and $Y$ with $X(x)=v$, $Y(x)=w$.
Let $\theta$ be a smooth $\R^k$-valued $1$-form on an open neighborhood $U$ of $x$ such that
$\Ker(\theta_x)=\mathcal D_x$ for all
$x\in U$. If $X$ and $Y$ are vector fields on an open neighborhood of $x$ then Cartan's formula
for exterior differentiation gives:
\[\dd\theta(X,Y)=X\big(\theta(Y)\big)-Y\big(\theta(X)\big)-\theta\big([X,Y]\big).\]
If $X$ and $Y$ are $\mathcal D$-horizontal then the equality above reduces to:
\[\dd\theta(X,Y)=-\theta\big([X,Y]\big).\]
The formula above implies that if $X'$, $Y'$ are $\mathcal D$-horizontal vector fields such that $X'(x)=X(x)$
and $Y'(x)=Y(x)$ then $\theta\big([X,Y]-[X',Y']\big)(x)=0$, i.e., $[X,Y](x)-[X',Y'](x)\in\mathcal D_x$.
Hence the Levi form is well-defined. Setting $X(x)=v$ and $Y(x)=w$ we obtain the following
formula:
\begin{equation}\label{eq:Levitheta}
\bar\theta_x\big(\Levi^{\mathcal D}_x(v,w)\big)=-\dd\theta(v,w),\quad v,w\in\mathcal D_x,
\end{equation}
where $\bar\theta_x:T_xE/\mathcal D_x\to\R^k$ denotes the linear map induced by $\theta_x$ in
the quotient space.

\begin{rem}
Clearly, by the classical Frobenius Theorem, $\mathcal D$ is integrable if and only if its Levi form
is identically zero. Moreover, the Levi form of $\mathcal D$ vanishes along any integral submanifold of $\mathcal D$.
\end{rem}

\begin{example}\label{exa:DGrafF}
Let $U$ be an open subset of $\R^n=\R^k\times\R^{n-k}$ and let
\[F:U\ni(x,y)\longmapsto F_{(x,y)}\in\Lin(\R^k,\R^{n-k})\]
be a smooth map. We consider the distribution $\mathcal D=\Gr(F)$ on $U$, i.e., $\mathcal D_{(x,y)}=\Gr\big(F_{(x,y)}\big)$,
for all $(x,y)\in U$. Given $X\in\R^k$, we define a $\mathcal D$-horizontal vector field $\widetilde X$ on $U$
by setting $\widetilde X_{(x,y)}=\big(X,F_{(x,y)}(X)\big)$, for all $(x,y)\in U$. Given $X,Y\in\R^k$ then:
\[[\widetilde X,\widetilde Y]=\big(0,\partial_xF(X,Y)+\partial_yF(F(X),Y)-
\partial_xF(Y,X)+\partial_yF(F(Y),X)\big).\]
If we identify $\mathcal D_{(x,y)}$ with $\R^k$ by the isomorphism $\big(X,F(X)\big)\mapsto X$ and
$\R^n/\mathcal D_{(x,y)}$ with $\R^{n-k}\cong\{0\}^k\times\R^{n-k}$ by the isomorphism
$(v,w)+\mathcal D_{(x,y)}\mapsto w-F(v)$ then the Levi form $\Levi^\mathcal D:\R^k\times\R^k\to\R^{n-k}$ is given
by:
\begin{equation}\label{eq:LeviGrF}
\Levi^\mathcal D(X,Y)=[\widetilde X,\widetilde Y].
\end{equation}
\end{example}

\begin{lem}\label{thm:switchHst}
Let $E$ be a smooth manifold, $\mathcal D$ be a smooth distribution on $E$ and let
\[U\ni(t,s)\longmapsto H(t,s)\in E\]
be a smooth map defined on an open subset $U\subset\R^2$. Let $I\subset\R$ be an interval
and let $s_0\in\R$ be such that $I\times\{s_0\}\subset U$ and $\Levi^{\mathcal D}_{H(t,s_0)}=0$
for all $t\in I$. Assume that $\frac{\partial H}{\partial t}(t,s)\in\mathcal D$ for all $(t,s)\in U$.
If $\frac{\partial H}{\partial s}(t_0,s_0)\in\mathcal D$ for some $t_0\in I$ then
$\frac{\partial H}{\partial s}(t,s_0)\in\mathcal D$ for all $t\in I$.
\end{lem}
\begin{proof}
The set:
\[I'=\big\{t\in I:\tfrac{\partial H}{\partial s}(t,s_0)\in\mathcal D\big\}\]
is obviously closed in $I$ because the map $I\ni t\mapsto\frac{\partial H}{\partial s}(t,s_0)\in TE$
is continuous and $\mathcal D$ is a closed subset of $TE$. Since $I$ is connected and $t_0\in I'$, the proof
will be complete once we show that $I'$ is open in $I$. Let $t_1\in I'$ be fixed. Let $\theta$
be an $\R^k$-valued smooth $1$-form defined in an open neighborhood $V$ of $H(t_1,s_0)$ in $E$
such that the linear map $\theta_x:T_xE\to\R^k$ is surjective and $\Ker(\theta_x)=\mathcal D_x$
for all $x\in V$. Choose a distribution $\mathcal D'$ on $V$ such that $T_xE=\mathcal D_x\oplus\mathcal D'_x$
for all $x\in V$. Then, for each $x\in V$, $\theta_x$ restricts to an isomorphism from $\mathcal D'_x$
onto $\R^k$. Let $J$ be a connected neighborhood of $t_1$ in $I$ such that
$H(t,s_0)\in V$ for all $t\in J$. We will show below that the map:
\begin{equation}\label{eq:themapHts0}
J\ni t\longmapsto\theta_{H(t,s_0)}\big(\tfrac{\partial H}{\partial s}(t,s_0)\big)\in\R^k
\end{equation}
is a solution of a homogeneous linear ODE; since
$\theta_{H(t_1,s_0)}\big(\frac{\partial H}{\partial
s}(t_1,s_0)\big)=0$, it will follow that
$\theta_{H(t,s_0)}\big(\frac{\partial H}{\partial
s}(t,s_0)\big)=0$ for all $t\in J$, i.e., $J\subset I'$.

We denote by $\frac\partial{\partial t}$ and $\frac\partial{\partial s}$ the canonical basis
of $\R^2$ and we apply Cartan's formula for exterior differentiation to the $1$-form
$H^*\theta$ obtaining:
\[\dd(H^*\theta)\big(\tfrac\partial{\partial t},\tfrac\partial{\partial s}\big)=
\tfrac\partial{\partial t}\Big((H^*\theta)\big(\tfrac\partial{\partial s}\big)\Big)-
\tfrac\partial{\partial s}\Big((H^*\theta)\big(\tfrac\partial{\partial t}\big)\Big)-
(H^*\theta)\big(\big[\tfrac\partial{\partial t},\tfrac\partial{\partial s}\big]\big).\]
Since $\dd(H^*\theta)=H^*(\dd\theta)$ and $\big[\frac\partial{\partial t},\frac\partial{\partial s}\big]=0$ we get:
\begin{align}\label{eq:thetaHts0}
\dd\theta_{H(t,s_0)}\big(\tfrac{\partial H}{\partial t}(t,s_0),\tfrac{\partial H}{\partial s}(t,s_0)\big)&=
\tfrac\partial{\partial t}\Big(\theta_{H(t,s_0)}\big(\tfrac{\partial H}{\partial s}(t,s_0)\big)\Big)\\
&-\tfrac\partial{\partial s}\big\vert_{s=s_0}\Big(\theta_{H(t,s)}\big(\tfrac{\partial H}{\partial t}(t,s)\big)\Big),
\quad t\in J.\notag
\end{align}
Observe that, since $\frac{\partial H}{\partial t}(t,s)$ is in $\mathcal D$, the last term
on the righthand side of \eqref{eq:thetaHts0} vanishes. We can write $\frac{\partial H}{\partial s}(t,s_0)=u_1(t)+u_2(t)$
with $u_1(t)\in\mathcal D$ and $u_2(t)\in\mathcal D'$. Since the Levi form of $\mathcal D$
vanishes at points of the form $H(t,s_0)$, equation \eqref{eq:Levitheta} implies that
$\dd\theta_{H(t,s_0)}(v,w)=0$ for all $v,w\in\mathcal D_{H(t,s_0)}$. We may thus replace
$\frac{\partial H}{\partial s}(t,s_0)$ by $u_2(t)$ in the lefthand side of \eqref{eq:thetaHts0}.
For $t\in J$ we consider the linear map $L(t):\R^k\to\R^k$ defined by:
\[L(t)\cdot z=\dd\theta_{H(t,s_0)}\big(\tfrac{\partial H}{\partial t}(t,s_0),\sigma_{H(t,s_0)}(z)\big),\quad z\in\R^k,\]
where, for $x\in V$, $\sigma_x:\R^k\to\mathcal D'_x$ denotes the inverse of the isomorphism
\[\theta_x\vert_{\mathcal D'_x}:\mathcal D'_x\longrightarrow\R^k.\]
Observe that:
\begin{align*}
\dd\theta_{H(t,s_0)}\big(\tfrac{\partial H}{\partial t}(t,s_0),\tfrac{\partial H}{\partial s}(t,s_0)\big)&=
\dd\theta_{H(t,s_0)}\big(\tfrac{\partial H}{\partial t}(t,s_0),u_2(t)\big)\\
&=L(t)\cdot\theta_{H(t,s_0)}\big(u_2(t)\big)\\
&=L(t)\cdot\theta_{H(t,s_0)}\big(\tfrac{\partial H}{\partial s}(t,s_0)\big).
\end{align*}
Equation \eqref{eq:thetaHts0} can now be rewritten as:
\[\tfrac{\partial}{\partial t}\Big(\theta_{H(t,s_0)}\big(\tfrac{\partial H}{\partial s}(t,s_0)\big)\Big)=
L(t)\cdot\theta_{H(t,s_0)}\big(\tfrac{\partial H}{\partial s}(t,s_0)\big),\quad t\in J.\]
Hence the map \eqref{eq:themapHts0} is a solution of a homogeneous linear ODE and we are done.
\end{proof}

\subsection{Horizontal distributions and horizontal liftings}

If $E$, $M$ are smooth manifolds and $\pi:E\to M$ is a smooth submersion then a smooth distribution $\mathcal D$
on $E$ is called {\em horizontal\/} with respect to $\pi$ if
\[T_xE=\Ker(\dd\pi_x)\oplus\mathcal D_x,\]
for all $x\in E$. Given a smooth horizontal distribution $\mathcal D$ on $E$ then a piecewise smooth
curve $\tilde\gamma:I\to E$ is called {\em horizontal\/} if $\tilde\gamma'(t)\in\mathcal D$ for all $t$ for which
$\tilde\gamma'(t)$ exists. Given a piecewise smooth curve $\gamma:I\to M$
then a {\em horizontal lifting\/} of $\gamma$ is a horizontal piecewise smooth curve
$\tilde\gamma:I\to E$ such that $\pi\circ\tilde\gamma=\gamma$.

By standard results of existence and uniqueness of solutions of ODE's it follows that given $t_0\in I$ and
$x_0\in\pi^{-1}\big(\gamma(t_0)\big)$ then there exists a unique maximal horizontal lifting $\tilde\gamma$ of $\gamma$
with $\tilde\gamma(t_0)=x_0$ defined in a subinterval of $I$ around $t_0$.

Let $\Lambda$ be a smooth manifold. By a {\em $\Lambda$-parametric family of curves\/} $\psi$ on $M$ we mean
a smooth map $\psi:Z\subset\R\times\Lambda\to M$ defined on an open subset $Z$ of $\R\times\Lambda$ such that the set:
\[I_\lambda=\big\{t\in\R:(t,\lambda)\in Z\big\}\subset\R\]
is an interval containing the origin, for all $\lambda\in\Lambda$. By a {\em local right inverse\/}
of $\psi$ we mean a locally defined smooth map $\alpha:V\subset M\to Z$ such that
$\psi\big(\alpha(m)\big)=m$, for all $m\in V$.

\begin{example}\label{exa:psigeodesic}
Let $M$ be a smooth manifold endowed with a connection $\nabla$. Given a point $x_0\in M$ we set $\Lambda=T_{x_0}M$ and we define
a $\Lambda$-parametric family of curves $\psi$ on $M$ by setting $\psi(t,\lambda)=\exp_{x_0}(t\lambda)$; the domain $Z\subset\R\times\Lambda$
of $\psi$ is the set of pairs $(t,\lambda)$ such that $t\lambda$ is in the domain of $\exp_{x_0}$. A local right inverse of $\psi$ is defined
as follows: let $V_0$ be an open neighborhood of the origin in $T_{x_0}M$ that is mapped diffeomorphically by $\exp_{x_0}$ onto an open
neighborhood $V$ of $x_0$ in $M$. We set:
\[\alpha(m)=\big(1,(\exp_{x_0}\vert_{V_0})^{-1}(m)\big),\]
for all $m\in V$. We remark that the same construction holds if one replaces the geodesic spray of a connection with
an arbitrary spray (see Section~\ref{sec:globalFrobenius}).
\end{example}

A {\em local section\/} of a smooth submersion $\pi:E\to M$ is a locally defined smooth map
$s:U\subset M\to E$ such that $\pi\circ s=\Id_U$. A local section $s$ is called {\em horizontal\/}
if the range of $\dd s(m)$ is $\mathcal D_{s(m)}$, for all $m\in U$.

\begin{lem}\label{thm:unicidadeFrobenius}
Let $s_1:U\to E$, $s_2:U\to E$ be local smooth horizontal sections of $E$ defined in an open connected subset
$U$ of $M$. If $s_1(x)=s_2(x)$ for some $x\in U$ then $s_1=s_2$.
\end{lem}
\begin{proof}
Given $y\in U$, there exists a piecewise smooth curve $\gamma:[a,b]\to U$ with $\gamma(a)=x$ and $\gamma(b)=y$.
Then $s_1\circ\gamma$ and $s_2\circ\gamma$ are both horizontal liftings of $\gamma$ starting at the same
point of $E$; hence $s_1\circ\gamma=s_2\circ\gamma$ and $s_1(y)=s_2(y)$.
\end{proof}

\begin{example}
Consider the distribution $\mathcal D=\Gr(F)$ on $U\subset\R^n$ defined in Example~\ref{exa:DGrafF}. Then
the first projection $\pi_1:U\to\R^k$ is a submersion and $\mathcal D$ is horizontal with respect to $\pi_1$.
A horizontal section $s:\R^k\supset V\to\R^n$ of $\pi_1$ is a map $s(x)=\big(x,f(x)\big)$ where
$f:V\to\R^{n-k}$ is a solution of the {\em total differential equation}:
\begin{equation}\label{eq:totalPDE}
\dd f(x)=F\big(x,f(x)\big),\quad x\in V.
\end{equation}
\end{example}

\subsection{The single leaf Frobenius theorem}

\begin{teo}[local single leaf Frobenius]\label{thm:Frobenius}
Let $E$, $M$ be smooth manifolds, $\pi:E\to M$ be a smooth submersion, $\mathcal D$ be a smooth horizontal distribution
on $E$ and $\psi:Z\subset\R\times\Lambda\to M$ be a $\Lambda$-parametric family of curves on $M$ with a local right
inverse $\alpha:V\subset M\to Z$. Let $\tilde\psi:Z\to E$ be a $\Lambda$-parametric family of curves on $E$ such that
$t\mapsto\tilde\psi(t,\lambda)$ is a horizontal lifting of $t\mapsto\psi(t,\lambda)$, for all $\lambda\in\Lambda$.
Assume that:
\begin{itemize}
\item[(a)] the Levi form of $\mathcal D$ vanishes on the range of $\tilde\psi$;
\item[(b)] $\partial_\lambda\tilde\psi(0,\lambda):T_\lambda\Lambda\to T_{\tilde\psi(0,\lambda)}E$ takes values
in $\mathcal D$ for all $\lambda\in\Lambda$.
\end{itemize}
Then $s=\tilde\psi\circ\alpha:V\to E$ is a local horizontal section of $\pi$.
\end{teo}
\begin{proof}
If $\left]-\varepsilon,\varepsilon\right[\ni s\mapsto\lambda(s)$ is an arbitrary smooth curve on $\Lambda$ then the map
\[H(t,s)=\tilde\psi\big(t,\lambda(s)\big)\]
satisfies the hypotheses of Lemma~\ref{thm:switchHst} with $t_0=0$ and $s_0=0$.
Thus:
\[\frac{\partial H}{\partial s}(t,0)=\frac{\partial\tilde\psi}{\partial\lambda}\big(t,\lambda(0)\big)\lambda'(0)\]
is in $\mathcal D$ for all $t\in I_{\lambda(0)}$. It follows that $\dd\tilde\psi_{(t,\lambda)}$ takes values
in $\mathcal D$, for all $(t,\lambda)\in Z$. Hence $\dd s(m)=\dd\tilde\psi\big(s(m)\big)\circ\dd\alpha(m)$
also takes values in $\mathcal D$, for all $m\in V$.
\end{proof}

\begin{rem}\label{thm:remtypical}
We observe that if the map $\lambda\mapsto\tilde\psi(0,\lambda)$ is constant then hypothesis (b) of Theorem~\ref{thm:Frobenius}
is automatically satisfied. Theorem~\ref{thm:Frobenius} is typically used as follows: one considers the $\Lambda$-parametric
family of curves $\psi$ explained in Example~\ref{exa:psigeodesic}, a fixed point $e_0\in\pi^{-1}(x_0)\subset E$ and for each $\lambda\in\Lambda$
one defines $t\mapsto\tilde\psi(t,\lambda)$ to be the horizontal lifting of $t\mapsto\psi(t,\lambda)$ with
$\tilde\psi(0,\lambda)=e_0$.
\end{rem}

\begin{example}\label{exa:FrobeniusGrF}
The single leaf Frobenius theorem can be used to prove the existence of solutions of the total
differential equation \eqref{eq:totalPDE} satisfying a initial condition $f(x_0)=y_0$ as follows.
Let $V$ be a star-shaped open neighborhood of $x_0$ in $\R^k$. Set $\Lambda=\R^k$; we define a
$\Lambda$-parametric family of curves $\psi:Z\subset\R\times\Lambda\to M$ on $M=\R^k$ by setting $\psi(t,\lambda)=
x_0+t\lambda$, where $Z$ is the set of pairs $(t,\lambda)$ with $x_0+t\lambda\in V$. A horizontal lifting
$t\mapsto\tilde\psi(t,\lambda)=\big(\psi(t,\lambda),\Psi(t,\lambda)\big)$ of the curve $t\mapsto\psi(t,\lambda)$
is a solution of the ODE:
\begin{equation}\label{eq:ODEpsiPsi}
\frac{\dd}{\dd t}\Psi(t,\lambda)=F_{\tilde\psi(t,\lambda)}(\lambda).
\end{equation}
We choose the solution $t\mapsto\tilde\psi(t,\lambda)$ of the ODE \eqref{eq:ODEpsiPsi} with initial condition
$\Psi(0,\lambda)=y_0$. We can assume that $V$ is small enough so that $\tilde\psi$ is well-defined on $Z$.
Hypothesis (b) of Theorem~\ref{thm:Frobenius} is then automatically satisfied and hypothesis (a) is equivalent to
the condition that \eqref{eq:LeviGrF} vanishes on the points of the form $\tilde\psi(t,\lambda)$, $(t,\lambda)\in Z$.
Under these circumstances, the thesis of Theorem~\ref{thm:Frobenius} guarantees that $f:V\ni x\mapsto\Psi(1,x-x_0)\in\R^{n-k}$
is a solution of the total differential equation \eqref{eq:totalPDE} with $f(x_0)=y_0$.
\end{example}

\subsection{The higher order single leaf Frobenius theorem}

Let $\mathcal D$ be a smooth distribution on a smooth manifold $E$. We denote by $\Gamma(TE)$ the set of all
smooth vector fields on $E$, by $\Gamma(\mathcal D)$ the subspace of $\Gamma(TE)$ consisting of
$\mathcal D$-horizontal vector fields and by $\Gamma^\infty(\mathcal D)$ the Lie subalgebra
of $\Gamma(TE)$ spanned by $\Gamma(\mathcal D)$. The Lie algebra $\Gamma^\infty(\mathcal D)$ can be alternatively
described as follows; we define recursively a sequence
\[\Gamma^0(\mathcal D)\subset\Gamma^1(\mathcal D)\subset\Gamma^2(\mathcal D)\subset\cdots\]
of subspaces of $\Gamma(TE)$ by setting $\Gamma^0(\mathcal D)=\Gamma(\mathcal D)$ and $\Gamma^{r+1}(\mathcal D)$ to be the
subspace of $\Gamma(TE)$ spanned by $\Gamma^r(\mathcal D)$ and by the brackets $[X,Y]$, with $X\in\Gamma^r(\mathcal D)$
and $Y\in\Gamma(\mathcal D)$. Then:
\[\Gamma^\infty(\mathcal D)=\bigcup_{r=0}^\infty\Gamma^r(\mathcal D).\]
Given $X\in\Gamma(TE)$ we denote by $\ad_X:\Gamma(TE)\to\Gamma(TE)$ the operator $\ad_X(Y)=[X,Y]$.

\begin{teo}\label{thm:higherorderFrobenius}
Let $E$ be a real-analytic manifold endowed with a real-analytic distribution $\mathcal D$. Given $e_0\in E$ then
there exists an integral submanifold of $\mathcal D$ passing through $e_0$ if and only if $X(e_0)\in\mathcal D_{e_0}$,
for all $X\in\Gamma^\infty(\mathcal D)$.
\end{teo}
\begin{proof}
If there exists an integral submanifold $S$ of $\mathcal D$ passing through $e_0$ then it follows immediately
by induction on $r$ that $X(S)\subset\mathcal D$, for all $X\in\Gamma^r(\mathcal D)$ and all $r\ge0$. Thus,
$X(e_0)\in\mathcal D_{e_0}$, for all $X\in\Gamma^\infty(\mathcal D)$. Conversely, assume that $X(e_0)\in\mathcal D_{e_0}$,
for all $X\in\Gamma^\infty(\mathcal D)$. By considering a convenient real-analytic local chart around $e_0$ we may assume
without loss of generality that $E=U$ is an open subset of $\R^n=\R^k\times\R^{n-k}$ and that $\mathcal D$ is of the
form $\Gr(F)$ (see Example~\ref{exa:DGrafF}). Write $e_0=(x_0,y_0)$; we will use the ideas explained in Example~\ref{exa:FrobeniusGrF}
to find a solution $f$ of the total differential equation \eqref{eq:totalPDE} with $f(x_0)=y_0$. Then $\Gr(f)$ is the required
integral submanifold of $\mathcal D$ passing through $e_0$. Observe that given $\lambda\in\R^k$ then
$t\mapsto\psi(t,\lambda)$ is an integral curve of the constant vector field $\lambda$ in $\R^k$ and thus
the horizontal lift $t\mapsto\tilde\psi(t,\lambda)$ is an integral curve of the vector field $\tilde\lambda=\big(\lambda,F(\lambda)\big)$
on $\R^n$ passing through $e_0$ at $t=0$. We now let $\lambda\in\R^k$, $X,Y\in\R^k$, be fixed and we define a map
$t\mapsto\phi(t)\in\R^{n-k}$ by setting:
\[\phi(t)=\Levi^\mathcal D_{\tilde\psi(t,\lambda)}(X,Y)=[\widetilde X,\widetilde Y]_{\psi(t,\lambda)}.\]
The proof will be completed once we show that $\phi$ is identically zero; since $\phi$ is real-analytic, it suffices
to proof that all derivatives of $\phi$ at $t=0$ vanish. Let us show by induction on $r$ that for all $r\ge0$
the $r$-th derivative of $\phi$ is given by:
\begin{equation}\label{eq:ufaphir}
\phi^{(r)}(t)=(\ad_{\tilde\lambda})^r[\widetilde X,\widetilde Y]+L^{(r)}\big(
(\ad_{\tilde\lambda})^i[\widetilde X,\widetilde Y];\ i=0,1,\ldots,r-1\big),
\end{equation}
where the righthand side is computed at the point $\tilde\psi(t,\lambda)$ and $L^{(r)}$ is a smooth map that
associates to each $(x,y)\in U\subset\R^n$ a linear map:
\[L^{(r)}_{(x,y)}:\bigoplus_r\R^{n-k}\longrightarrow\R^{n-k}.\]
From equality \eqref{eq:ufaphir} the conclusion will follow; namely, for all $i$,
$(\ad_{\tilde\lambda})^i[\widetilde X,\widetilde Y]$ is in $\{0\}^k\times\R^{n-k}$ and
since $\big((\ad_{\tilde\lambda})^i[\widetilde X,\widetilde Y]\big)_{e_0}\in\mathcal D_{e_0}$,
we get $\big((\ad_{\tilde\lambda})^i[\widetilde X,\widetilde Y]\big)_{e_0}=0$. Hence $\phi^{(r)}(0)=0$, for all $r\ge0$.
To prove \eqref{eq:ufaphir} simply differentiate both sides with respect to $t$, observing that:
\[\frac{\dd}{\dd t}(\ad_{\tilde\lambda})^i[\widetilde X,\widetilde Y]=\dd\big((\ad_{\tilde\lambda})^i[\widetilde X,\widetilde Y]\big)
\cdot\tilde\lambda=(\ad_{\tilde\lambda})^{i+1}[\widetilde X,\widetilde Y]+
\dd\tilde\lambda\big((\ad_{\tilde\lambda})^i[\widetilde X,\widetilde Y]\big).\]
\end{proof}

\begin{rem}\label{thm:highorderref}
Clearly, the hypotheses of Theorem~\ref{thm:higherorderFrobenius} are local, i.e., if $U$ is an open neighborhood
of $e_0$ in $E$ then $X(e_0)\in\mathcal D_{e_0}$ for all $X\in\Gamma^\infty(\mathcal D\vert_U)$ if and only if
$X(e_0)\in\mathcal D_{e_0}$ for all $X\in\Gamma^\infty(\mathcal D)$. Replacing $E$ with an open neighborhood of $e_0$,
we may assume that $\mathcal D$ admits a global referential $X_1$, \dots, $X_k$. It is easy to see
that $\Gamma^r(\mathcal D)$ is the $C^\infty(E)$-module spanned by $X_1$, \dots, $X_k$ and by the iterated brackets:
\begin{equation}\label{eq:rorderbrack}
[X_{i_1},[X_{i_2},\ldots,[X_{i_s},X_{i_{s+1}}]\cdots]],\quad i_1,\ldots,i_s=1,\ldots,k,\ s=1,\ldots,r.
\end{equation}
Thus, in order to check the hypotheses of Theorem~\ref{thm:higherorderFrobenius},
it suffices to verify if the brackets in \eqref{eq:rorderbrack}
evaluated at $e_0$ are in $\mathcal D_{e_0}$, for all $s\ge1$.
\end{rem}

\end{section}

\begin{section}{The global ``single leaf Frobenius Theorem''}
\label{sec:globalFrobenius}

\subsection{Sprays on manifolds}\label{sub:sprays}
Let $M$ be a smooth manifold and let $\pi:TM\to M$ the canonical projection of
its tangent bundle. Denote by $\dd\pi:TTM\to TM$ the
differential of $\pi$; we denote by $\bar\pi:TTM\to TM$ the natural projection of $TTM=T(TM)$. For each
$a\in\R$ we denote by $\mathfrak m_a:TM\to TM$ the operator of multiplication by $a$.
\begin{defin}\label{thm:defspray}
A {\em spray\/} on $M$ is a smooth vector field $\mathcal S:TM\to TTM$ on the manifold $TM$ satisfying the following two
conditions:
\begin{itemize}
\item[(i)] $\dd\pi\circ\mathcal S=\bar\pi\circ\mathcal S$;
\item[(ii)] for all $a\in\R$, $a\,\dd\mathfrak m_a\circ\mathcal S=\mathcal S\circ\mathfrak m_a$, i.e.,
$a\,\dd\mathfrak m_a(v)\mathcal S(v)=\mathcal S(av)$, for all $v\in TM$.
\end{itemize}
\end{defin}

\begin{rem}\label{thm:remzerosection}
Notice that property (b) on Definition~\ref{thm:defspray} implies that a spray vanishes on the zero section
of $TM$. In particular, the integral curves of $\mathcal S$ passing through the zero section are constant.
\end{rem}

\begin{lem}
Let $\mathcal S:TM\to TTM$ be a smooth vector field on $TM$. Then $\mathcal S$ is a spray on $M$ if and only if
the following conditions are satisfied:
\begin{itemize}
\item[(a)] for every integral curve $\lambda:I\to TM$ of $\mathcal S$, we have $\lambda=\gamma'$, where $\gamma=\pi\circ\lambda$;
\item[(b)] if $\lambda=\gamma':I\to TM$ is an integral curve of $\mathcal S$ then
\[I\ni t\longmapsto\frac{\dd}{\dd t}\gamma(at)\in TM\]
is an integral curve of $\mathcal S$, for all $a\in\R$.
\end{itemize}
\end{lem}

\begin{defin}
A curve $\gamma:I\to M$ is called a (maximal) solution of $\mathcal S$ if $\gamma':I\to TM$ is a (maximal) integral
curve of the vector field $\mathcal S$.
\end{defin}
Obviously for every $x\in M$, $v\in T_xM$ there exists a unique maximal solution $\gamma$ of $\mathcal S$ with
$\gamma(0)=x$ and $\gamma'(0)=v$.

\begin{example}[geodesic spray]\label{exa:geodesicspray}
If $\nabla$ is a connection on $M$ then one can define a spray $\mathcal S$ on $M$ by taking $\mathcal S(v)$
to be the unique horizontal vector on $T_vTM$ such that $\dd\pi_v\big(\mathcal S(v)\big)=v$, for all $v\in TM$.
The integral curves of $\mathcal S$ are the curves $\gamma'$, with $\gamma:I\to M$ a geodesic of $\nabla$.
\end{example}

\begin{example}[one-parameter subgroup spray]
Let $G$ be a Lie group and denote by $\mathfrak g$ its Lie algebra.
Using left (resp., right) translations, one can identify the tangent bundle
$TG$ with the product $G\times\mathfrak g$, so that
\[T(TG)\cong T(G\times\mathfrak g)\cong(TG)\times(T\mathfrak g)\cong(G\times\mathfrak g)\times(\mathfrak g\times\mathfrak g).\]
The vector field on $TG$ given by $\mathcal S(g,X)=(g,X,X,0)$, $g\in G$, $X\in\mathfrak g$,
is a spray in $G$, whose solutions are left (resp., right) translations of
one-parameter subgroups of $G$. The spray $\mathcal S$ is the geodesic spray of the connection whose Christoffel symbols vanish
on a left (resp., right) invariant frame.
\end{example}

Let $\mathcal S$ be a fixed spray on $M$ and denote by
\[F:\Dom(F)\subset\R\times TM\longrightarrow TM\]
its maximal flow. The {\em exponential map\/} associated to $\mathcal S$ is the map:
\[\exp(v)=\pi\big(F(1,v)\big)\in M,\]
defined on the set:
\[\Dom(\exp)=\big\{v\in TM:(1,v)\in\Dom(F)\big\}.\]
Since $\Dom(F)$ is open in $\R\times TM$, $\Dom(\exp)$ is open in $TM$; moreover, by Remark~\ref{thm:remzerosection}
the zero section of $TM$ is contained in $\Dom(\exp)$. In particular, for each $x\in M$, the intersection of $\Dom(\exp)$
with $T_xM$ is an open neighborhood of the origin.

\begin{lem}
For all $t,s\in\R$, $v\in TM$, $(t,sv)\in\R\times TM$ is in $\Dom(F)$ if and only if $(ts,v)\in\R\times TM$ is in $\Dom(F)$;
moreover, $F(t,sv)=sF(ts,v)$.
\end{lem}

\begin{cor}
For all $s\in\R$, $v\in TM$, $(s,v)\in\R\times TM$ is in $\Dom(F)$ if and only if $sv$ is in $\Dom(\exp)$;
moreover, $\pi\big(F(s,v)\big)=\exp(sv)$.
\end{cor}

\begin{cor}\label{thm:starshaped}
Given $x\in M$, $v\in T_xM$ then the set $\big\{s\in\R:sv\in\Dom(\exp)\big\}$ is an open interval containing the origin;
the map $\gamma(s)=\exp(sv)$ defined on such open interval is the maximal solution of $\mathcal S$
with $\gamma(0)=x$, $\gamma'(0)=v$.
\end{cor}

For each $x\in M$ let us denote by $\exp_x$ the restriction of $\exp$ to $\Dom(\exp)\cap T_xM$.
It follows from Corollary~\ref{thm:starshaped} that the domain of $\exp_x$ is a star-shaped open neighborhood
of the origin in $T_xM$; moreover, $\dd\exp_x(0)$ is the identity map of $T_xM$.

\begin{defin}
A {\em normal neighborhood\/} of a point $x\in M$ is an open neighborhood $V\subset M$ of $x$ such that
there exists a star-shaped open neighborhood $U$ of the origin in $T_xM$ such that $\exp_x\vert_U:U\to V$
is a diffeomorphism. An open subset $V$ of $M$ is called {\em normal\/}\footnote{%
Observe that, according to this definition, a normal open subset of $M$ containing a point $x\in M$ is not
necessarily a normal neighborhood of $x$!} if every $x\in M$ has a normal
neighborhood containing $V$.
\end{defin}

It follows from the inverse function theorem that every point of $M$ has a normal neighborhood. Moreover, we have
the following:
\begin{prop}\label{thm:propopennormal}
Every point of $M$ is contained in some normal open subset of $M$.
\end{prop}
\begin{proof}
Consider the map $\phi:\Dom(\exp)\subset TM\to M\times M$ given by $\phi(v)=\big(\exp(v),\pi(v)\big)$.
Given $x\in M$ and denote by $0_x\in TM$ the origin of $T_xM$. We identify $T_{0_x}TM$ with $T_xM\oplus T_xM$,
where the first summand corresponds to the tangent space of the zero section of $TM$ and the second summand
corresponds to the tangent space to the fiber of $TM$ containing $0_x$. The differential of $\phi$
at $0_x$ is easily computed as:
\[\dd\phi_{0_x}(v,w)=(v+w,v),\quad v,w\in T_xM.\]
It follows from the inverse function theorem that $\phi$ carries an open neighborhood $\mathcal U$
of $0_x$ in $TM$ diffeomorphically onto an open neighborhood of $(x,x)$ in $M\times M$.
We can choose $\mathcal U$ such that $\mathcal U\cap T_yM$ is a star-shaped open neighborhood of the origin
of $T_yM$, for all $y\in\pi(\mathcal U)$. Let $V$ be an open neighborhood of $x$ in $M$ such that
$V\times V\subset\phi(\mathcal U)$. We claim that $V$ is a normal open subset of $M$.
Let $y\in V$ be fixed. Clearly $V\subset\pi(\mathcal U)$, so that $\mathcal U\cap T_yM$ is a star-shaped
open neighborhood of the origin of $T_yM$; thus $\exp(\mathcal U\cap T_yM)$ is a normal neighborhood of $y$.
Moreover, given $z\in V$ then $(z,y)\in V\times V$, so that there exists $v\in\mathcal U$ with $\phi(v)=(z,y)$;
then $v\in\mathcal U\cap T_yM$ and hence $z\in\exp(\mathcal U\cap T_yM)$.
\end{proof}

\begin{defin}
A {\em piecewise solution\/} of a spray $\mathcal S$ is a curve $\gamma:[a,b]\to M$ for which there exists
a partition $a=t_0<t_1<\cdots<t_k=b$ of $[a,b]$ such that $\gamma\vert_{[t_i,t_{i+1}]}$ is a solution
of $\mathcal S$ for all $i$.
\end{defin}

\subsection{The global single leaf Frobenius theorem}

\begin{teo}[global single leaf Frobenius]\label{thm:globalFrobenius}
Let $E$, $M$ be smooth manifolds, $\pi:E\to M$ be a smooth submersion and $\mathcal D$ be a smooth horizontal distribution
on $E$. Let $x_0\in M$, $e_0\in\pi^{-1}(x_0)\subset E$ be given and let $\mathcal S$ be a fixed spray on $M$.
Assume that:
\begin{itemize}
\item[(a)] every piecewise solution $\gamma:[a,b]\to M$ of $\mathcal S$ with $\gamma(a)=x_0$ admits a horizontal lifting
$\tilde\gamma:[a,b]\to E$ with $\tilde\gamma(a)=e_0$;
\item[(b)] if $\tilde\gamma:[a,b]\to E$ is the horizontal lifting of a piecewise solution $\gamma:[a,b]\to M$ of $\mathcal S$
with $\tilde\gamma(a)=e_0$ then the Levi form of $\mathcal D$ vanishes at the point $\tilde\gamma(b)\in E$;
\item[(c)] $M$ is (connected and) simply-connected.
\end{itemize}
Then there exists a unique global smooth horizontal section $s$ of $E$ with $s(x_0)=e_0$.
\end{teo}
\begin{proof}
Uniqueness follows directly from Lemma~\ref{thm:unicidadeFrobenius}. For the existence,
we use the globalization theory explained in Appendix~\ref{sec:globalization}.

Let $E'$ denote the subset of $E$ consisting of the points of the form $\tilde\gamma(b)$, where $\tilde\gamma(a)=e_0$
and $\tilde\gamma:[a,b]\to E$ is the horizontal lifting of some piecewise solution $\gamma:[a,b]\to M$ of $\mathcal S$ with $\gamma(a)=x_0$.
We define a pre-sheaf $\mathfrak P$ over $M$ as follows: for each open subset $U$ of $M$, $\mathfrak P(U)$ is the set
of all smooth horizontal sections $s:U\to E$ with $s(U)\subset E'$. Given open subsets $U,V\subset M$ with $V\subset U$
then $\mathfrak P_{U,V}$ is given by $\mathfrak P_{U,V}(s)=s\vert_V$, for all $s\in\mathfrak P(U)$.
The existence of a global smooth horizontal section of $E$ is equivalent to $\mathfrak P(M)\ne\emptyset$.
We will use Proposition~\ref{thm:niceglobalization}. Using Theorem~\ref{thm:Frobenius} (recall Remark~\ref{thm:remtypical})
we get a smooth horizontal section $s:U\to E$ defined in an open neighborhood $U$ of $x_0$; it is clear by the construction
of $s$ that $s(U)\subset E'$. Thus the pre-sheaf $\mathfrak P$ is nontrivial. The localization property (Definition~\ref{thm:localizationproperty})
for $\mathfrak P$
is trivial and the uniqueness property (Definition~\ref{thm:uniquenessprop})
for $\mathfrak P$ follows directly from Lemma~\ref{thm:unicidadeFrobenius}.
To conclude the proof, we show that $\mathfrak P$ has the extension property (Definition~\ref{thm:uniquenessprop}).
We shall prove that every normal open subset of $M$ has the extension property for $\mathfrak P$ (recall Proposition~\ref{thm:propopennormal}).
Let $U$ be an open normal subset of $M$, $V$ be a nonempty open connected subset of $U$ and $s\in\mathfrak P(V)$
be a smooth horizontal section of $E$ with $s(V)\subset E'$. Let $x_1\in V$ be fixed. Since $s(x_1)\in E'$, there
exists a piecewise solution $\gamma:[a,b]\to M$ of $\mathcal S$ with $\gamma(a)=x_0$ and a horizontal lifting
$\tilde\gamma:[a,b]\to E$ of $\gamma$ with $\tilde\gamma(a)=e_0$ and $\tilde\gamma(b)=s(x_1)$. Let $W$
be a normal neighborhood of $x_1$ containing $U$ and $W_0$ be a star-shaped open neighborhood of the origin
in $T_{x_1}M$ such that $\exp_{x_1}:W_0\to W$ is a diffeomorphism. For each $x\in W$ let $v\in W_0$ be such that
$\exp_{x_1}(v)=x$; we claim that $\mu_x:[0,1]\ni t\mapsto\exp_{x_1}(tv)\in M$ has a horizontal lifting $\tilde\mu:[0,1]\to E$
starting at $s(x_1)$ and that the Levi form of $\mathcal D$ vanishes along the image of $\tilde\mu$. Namely,
the concatenation $\gamma\cdot\mu$ of $\gamma$ with $\mu$ is a piecewise solution of $\mathcal S$ starting
at $x_0$; by hypothesis~(a), $\gamma\cdot\mu$ has a horizontal lifting starting at $e_0$. Such horizontal lifting
is of the form $\tilde\gamma\cdot\tilde\mu$, where $\tilde\mu$ is a horizontal lifting of $\mu$ starting at $s(x_1)$;
moreover, hypothesis~(b) implies that the Levi form of $\mathcal D$ vanishes along $\tilde\gamma\cdot\tilde\mu$.
Observe that the image of $\tilde\mu$ is contained in $E'$. We can now apply Theorem~\ref{thm:Frobenius} to obtain
a smooth horizontal section $\bar s:W\to E$ with $\bar s(x_1)=s(x_1)$. Thus, by Lemma~\ref{thm:unicidadeFrobenius}
and the connectedness of $V$, $\bar s\vert_V=s$ and hence $\bar s\vert_U\in\mathfrak P(U)$ is an extension of $s$
to $U$.
\end{proof}

\begin{prop}\label{thm:Frobeniusrealanalytic}
Let $E$, $M$ be real-analytic manifolds, $\pi:E\to M$ be a real-analytic submersion and $\mathcal D$ be a real-analytic
horizontal distribution on $E$. Assume that:
\begin{itemize}
\item[(a)] $M$ is (connected and) simply-connected;
\item[(b)] given a real analytic curve $\gamma:I\to M$, $t_0\in I$ and $e_0\in\pi^{-1}\big(\gamma(t_0)\big)$ then
there exists a horizontal lifting $\tilde\gamma:I\to E$ of $\gamma$ with $\tilde\gamma(t_0)=e_0$.
\end{itemize}
Then any local horizontal section $s:U\to E$ of $\pi$ defined on a nonempty connected open subset $U$ of $M$
extends to a global horizontal section of $\pi$. In particular, if $\mathcal D$ satisfies the hypothesis of
Theorem~\ref{thm:higherorderFrobenius} at some point $e_0$ of $E$, assumptions (a) and (b) imply that $\pi$ admits
a global horizontal section.
\end{prop}
\begin{proof}
We use again the globalization theory explained in Appendix~\ref{sec:globalization}.
We define a pre-sheaf $\mathfrak P$ over $M$ as follows: for each open subset $U$ of $M$, $\mathfrak P(U)$ is the set
of all smooth horizontal sections $s:U\to E$; given open subsets $U,V\subset M$ with $V\subset U$
then $\mathfrak P_{U,V}$ is given by $\mathfrak P_{U,V}(s)=s\vert_V$, for all $s\in\mathfrak P(U)$. By Proposition~\ref{thm:niceglobalization}
it suffices to show that $\mathfrak P$ has the localization property, the uniqueness property and the extension property.
The localization property is trivial and the uniqueness property follows from Lemma~\ref{thm:unicidadeFrobenius}.
As to the extension property, it can be proved as follows. Let $x_0\in M$ be fixed and let $\varphi:U\to\B_0(r)$
be a real-analytic chart defined on an open neighborhood $U$ of $x_0$, taking values in the open ball
$\B_0(r)\subset\R^n$ of radius $r$ centered at the origin and $\varphi(x_0)=0$. We will show that $V=\varphi^{-1}\big(\B_0(r/3)\big)$
is an open neighborhood of $x_0$ having the extension property for $\mathfrak P$. To this aim, let $W$ be a nonempty
connected open subset of $V$ and let $s\in\mathfrak P(W)$ be a local horizontal section defined on $W$.
Choose $x_1\in W$. Set $\Lambda=\B_0\big(\frac23r\big)$ and let $\psi:Z\subset\R\times\Lambda\to M$ be the one-parameter
family of curves defined by $\psi(t,\lambda)=\varphi^{-1}\big(\varphi(x_1)+t\lambda\big)$, where $Z$ is the
set of pairs $(t,\lambda)\in\R\times\Lambda$ with $\varphi(x_1)+t\lambda\in\B_0(r)$. We define a local right inverse
\[\alpha:\varphi^{-1}\big(\B_{\varphi(x_1)}(\tfrac23r)\big)\longrightarrow Z\subset\R\times\Lambda\]
of $\psi$ by setting $\alpha(x)=\big(1,\varphi(x)-\varphi(x_1)\big)$. By assumption~(b), for each $\lambda\in\Lambda$,
the curve $t\mapsto\psi(t,\lambda)$ has a horizontal lifting $t\mapsto\tilde\psi(t,\lambda)\in E$ with
$\tilde\psi(0,\lambda)=s(x_1)$. Notice that, by the uniqueness of the horizontal lifting of a curve, we have
$\tilde\psi(t,\lambda)=s\big(\psi(t,\lambda)\big)$ for small $t$.
Since $s$ is a horizontal section of $\pi$, its image is an integral submanifold
of $\mathcal D$ and thus the Levi form $\Levi^\mathcal D$ vanishes along the image of $s$. Thus $\Levi^\mathcal D$
vanishes at the point $\tilde\psi(t,\lambda)$ for small $t$; hence, since $t\mapsto\Levi^\mathcal D\big(\tilde\psi(t,\lambda)\big)$
is real-analytic, $\Levi^\mathcal D$ must vanish along the entire curve $t\mapsto\tilde\psi(t,\lambda)$.
By Theorem~\ref{thm:Frobenius}, $\tilde\psi\circ\alpha$ is a horizontal section of $E$ with $(\tilde\psi\circ\alpha)(x_1)=s(x_1)$;
since the domain of $\alpha$ clearly contains $V$, Lemma~\ref{thm:unicidadeFrobenius} implies that $\tilde\psi\circ\alpha$
extends $s$ to (an open set containing) $V$. This proves the extension property of $\mathfrak P$ and concludes the proof.
\end{proof}

\end{section}

\begin{section}{Levi--Civita connections}

\begin{subsection}{Levi form of the horizontal distribution of a connection}

Let $\pi:E\to M$ be a smooth vector bundle over a smooth manifold
$M$ and let $\nabla$ be a connection on $E$; for $m\in M$ we denote
by $E_m=\pi^{-1}(m)$ the fiber of $E$ over $m$. We denote by $R_m:T_mM\times T_mM\times E_m\to E_m$
the {\em curvature tensor\/} of $\nabla$ defined by:
\[R(X,Y)\xi=\nabla_X\nabla_Y\xi-\nabla_Y\nabla_X\xi-\nabla_{[X,Y]}\xi,\]
for all smooth vector fields $X$, $Y$ in $M$ and every smooth section $\xi$ of $E$.

Recall that there
exists a unique distribution $\mathcal D$ on the manifold $E$ that
is horizontal with respect to $\pi$ and has the following property:
if $\gamma:I\to M$ is a smooth curve on $M$ then a curve
$\tilde\gamma:I\to E$ is a horizontal lifting of $\gamma$ with
respect to $\mathcal D$ if and only if $\tilde\gamma$ is a
$\nabla$-parallel section of $E$ along $\gamma$. We call $\mathcal
D$ the {\em horizontal distribution\/} of $\nabla$. Given $m\in M$
and $\xi\in E_m$ then the quotient $T_\xi E/\mathcal D_\xi$ can be
identified with $T_\xi(E_m)=\Ker(\dd\pi_\xi)$; moreover, since $E_m$
is a vector space, we identify $T_\xi(E_m)$ with $E_m$. We also
identify $\mathcal D_\xi$ with $T_mM$ using $\dd\pi_\xi$. The Levi
form of $\mathcal D$ at a point $\xi\in E$ can thus be seen as a
bilinear map:
\[\Levi^{\mathcal D}_\xi:T_mM\times T_mM\longrightarrow E_m.\]

\begin{lem}\label{thm:Levihorizontal}
The Levi form of the horizontal distribution $\mathcal D$ of a connection $\nabla$ is given by:
\[\Levi^{\mathcal D}_\xi(v,w)=-R_m(v,w)\xi,\]
for all $m\in M$, $\xi\in E_m$.
\end{lem}
\begin{proof}
Given a smooth vector field $X$ on $M$ we denote by $X^\hor$ the {\em horizontal lift\/} of $X$ which is the
unique horizontal vector field on $E$ such that $\dd\pi_\xi(X^\hor(\xi)\big)=X\big(\pi(\xi)\big)$, for all $\xi\in E$.
Given smooth vector fields $X$, $Y$ on $M$, we have to show that vertical component of $[X^\hor,Y^\hor]$ at
a point $\xi\in E$ is equal to $-R(X,Y)\xi$. Note that the horizontal component of $[X^\hor,Y^\hor]$ is
$[X,Y]^\hor$, since $X^\hor$ and $Y^\hor$ are $\pi$-related respectively with $X$ and $Y$. Thus, the proof will
be concluded once we show that:
\[\alpha\big([X^\hor,Y^\hor]-[X,Y]^\hor\big)=-\alpha\big(R(X,Y)\xi\big),\]
for every smooth section $\alpha$ of the dual bundle $E^*$. Given one such section $\alpha$, we denote by $f_\alpha:E\to\R$
the smooth map defined by:
\[f_\alpha(\xi)=\alpha(\xi).\]
We claim that:
\[X^\hor(f_\alpha)=f_{\lower2pt\hbox{$\scriptstyle\nabla^*_{\!X}\alpha$}},\]
where $\nabla^*$ denotes the connection of $E^*$.
Namely, let $\gamma:\left]-\varepsilon,\varepsilon\right[\to M$ be an integral curve of $X$ and let $t\mapsto\xi(t)$
be a parallel section of $E$ along $\gamma$, so that $\xi$ is an integral curve of $X^\hor$; then:
\begin{multline*}
X^\hor(f_\alpha)=\tfrac\dd{\dd t}\big\vert_{t=0}f_\alpha\big(\xi(t)\big)=\tfrac\dd{\dd t}\big\vert_{t=0}
\alpha_{\gamma(t)}\big(\xi(t)\big)\\
=\big(\nabla^*_{\gamma'(t)}\alpha\big)\xi(t)\big\vert_{t=0}=(\nabla^*_X\alpha)\xi,
\end{multline*}
which proves the claim. Observe also that if $v\in TE$ is a vertical vector then $v(f_\alpha)=\alpha(v)$; therefore:
\begin{equation}\label{eq:bothsides}
\alpha\big([X^\hor,Y^\hor]-[X,Y]^\hor\big)=\big([X^\hor,Y^\hor]-[X,Y]^\hor\big)(f_\alpha)
=f_{R^*(X,Y)\alpha},
\end{equation}
where $R^*$ denotes the curvature tensor of $\nabla^*$. A simple computation shows that:
\[R^*(X,Y)\alpha=-\alpha\circ R(X,Y).\]
The conclusion follows from \eqref{eq:bothsides} by evaluating both sides at the point $\xi$.
\end{proof}

\begin{cor}\label{thm:corparalelo}
Let $\pi:E\to M$ be a smooth vector bundle endowed with a connection $\nabla$, let
$\psi:Z\subset\R\times\Lambda\to M$ be a $\Lambda$-parametric
family of curves on $M$ with a local right inverse $\alpha:V\subset M\to Z$ and let $\tilde\psi:Z\to E$ be a smooth
section of $E$ along $\psi$ such that $t\mapsto\tilde\psi(t,\lambda)$ is parallel for all $\lambda\in\Lambda$ and such
that $\lambda\mapsto\tilde\psi(0,\lambda)$ is also parallel.
If
\[R_{\psi(t,\lambda)}(v,w)\tilde\psi(t,\lambda)=0,\]
for all $v,w\in T_{\psi(t,\lambda)}M$ and all
$(t,\lambda)\in Z$ then $s=\tilde\psi\circ\alpha$ is a parallel local section of $E$.
\end{cor}
\begin{proof}
Follows readily from Theorem~\ref{thm:Frobenius} and Lemma~\ref{thm:Levihorizontal}.
\end{proof}

\begin{cor}\label{thm:corparaleloglobal}
Let $\pi:E\to M$ be a smooth vector bundle endowed with a connection $\nabla$.
Let $x_0\in M$, $e_0\in\pi^{-1}(x_0)\subset E$ be given and let $\mathcal S$ be a fixed spray on $M$.
Assume that:
\begin{itemize}
\item[(a)] if $\gamma:[a,b]\to M$ is a piecewise solution of $\mathcal S$ with $\gamma(a)=x_0$ and
$\tilde\gamma:[a,b]\to E$ is a parallel section of $E$ along $\gamma$ with $\tilde\gamma(a)=e_0$ then
$R_{\gamma(b)}(v,w)\tilde\gamma(b)=0$, for all $v,w\in T_{\gamma(b)}M$;
\item[(b)] $M$ is (connected and) simply-connected.
\end{itemize}
Then there exists a unique global smooth parallel section $s$ of $E$ with $s(x_0)=e_0$.
\end{cor}
\begin{proof}
Follows directly from Lemma~\ref{thm:Levihorizontal} and Theorem~\ref{thm:globalFrobenius}.
\end{proof}

\begin{cor}\label{thm:correalanalyticparallel}
Let $\pi:E\to M$ be a real-analytic vector bundle endowed with a real-analytic connection $\nabla$. Assume that
$M$ is (connected and) simply-connected.
Then any local parallel section $s:U\to E$ of $E$ defined on a nonempty connected open subset $U$ of $M$
extends to a global parallel section of $E$.
\end{cor}
\begin{proof}
It follows from Lemma~\ref{thm:Levihorizontal} and Proposition~\ref{thm:Frobeniusrealanalytic}.
\end{proof}

\begin{prop}\label{thm:prophighorderhorizontal}
Let $\pi:E\to M$ be a real-analytic vector bundle endowed with a real-analytic connection $\nabla$.
Given $x\in M$, $e\in\pi^{-1}(x)$, assume that:
\begin{equation}\label{eq:condnablakR}
(\nabla^kR)(v_1,v_2,\ldots,v_{k+2})e=0,
\end{equation}
for all $v_1,\ldots,v_{k+2}\in T_xM$ and all $k\ge0$.
Then there exists a parallel section $s$ of $E$ defined in an open neighborhood of $x$ in $M$ with $s(x)=e$;
in particular, by Corollary~\ref{thm:correalanalyticparallel}, if $M$ is (connected and) simply-connected then
there exists a global parallel section $s$ of $E$ with $s(x)=e$.
\end{prop}
\begin{proof}
Given a smooth vector field $X$ on $M$, we denote by $\widehat X$ the unique horizontal vector field on $E$
that is $\pi$-related with $X$.
We show that condition \eqref{eq:condnablakR} is equivalent to the condition that all iterated
brackets of vector fields $\widehat X$ are horizontal at the point $e$. The conclusion will then follow
from Theorem~\ref{thm:higherorderFrobenius}. First, let us compute the bracket $[\widehat X,\widehat Y]$.
Since $\widehat X$ and $\widehat Y$ are $\pi$-related respectively with $X$ and $Y$, it follows that the horizontal
component of $[\widehat X,\widehat Y]$ is $[X,Y]$; its vertical component is computed in Lemma~\ref{thm:Levihorizontal}.
Thus:
\begin{equation}\label{eq:brachatXhatY}
[\widehat X,\widehat Y]_e=\big([X,Y]_x,-R(X,Y)e\big),
\end{equation}
where we write tangent vectors to $E$ as pairs consisting of a horizontal component and a vertical component.
Given a smooth section $L$ of the vector bundle $\Lin(E)$, we denote by $\widetilde L$ the vertical vector field
on $E$ defined by $\widetilde L(e)=\big(0,L(e)\big)$. Given a smooth vector field $Z$ on $M$, let us compute
the bracket $[\widehat Z,\widetilde L]$. Since $\widehat Z$ is $\pi$-related with $Z$ and $\widetilde L$ is $\pi$-related
with zero, it follows that $[\widehat Z,\widetilde L]$ is vertical. Given a smooth section $\alpha$ of $E^*$, we consider
the map $f_\alpha:E\to\R$ defined by $f_\alpha(e)=\alpha(e)$ and we compute as follows:
\begin{gather*}
\widetilde L(f_\alpha)(e)=\alpha\big(L(e)\big)=f_{\alpha\circ L}(e),\\
\widehat Z(f_\alpha)(e)=\frac{\dd}{\dd t}f_\alpha\big(e(t)\big)=\frac{\dd}{\dd t}\,\alpha\big(e(t)\big)
=(\nabla_Z\alpha)(e)=f_{\nabla_Z\alpha}(e),
\end{gather*}
where $t\mapsto e(t)$ is an integral curve of $\widehat Z$, i.e., a parallel section of $E$ along an integral curve of
$Z$. Then:
\begin{multline*}[\widehat Z,\widetilde L](f_\alpha)=\widehat Z\big(\widetilde L(f_\alpha)\big)-
\widetilde L\big(\widehat Z(f_\alpha)\big)=f_{\nabla_Z(\alpha\circ L)}-f_{(\nabla_Z\alpha)\circ L}
=f_{\alpha\circ\nabla_ZL}\\
=\widetilde{\nabla_ZL}(f_\alpha),
\end{multline*}
so that:
\begin{equation}\label{eq:brackhatZtildeL}
[\widehat Z,\widetilde L]=\widetilde{\nabla_ZL}.
\end{equation}
Notice that \eqref{eq:brachatXhatY} says that $[\widehat X,\widehat Y]$ is given by:
\[[\widehat X,\widehat Y]=\widehat{[X,Y]}-\widetilde L,\]
where $L(e)=R(X,Y)e$. Using the equality above and \eqref{eq:brackhatZtildeL} it can be easily proved by induction
that:
\[[\widehat Z_1,[\widehat Z_2,\ldots[\widehat Z_k,[\widehat X,\widehat Y]]\cdots]]=
[\widehat Z_1,[\widehat Z_2,\ldots[\widehat Z_k,\widehat{[X,Y]}]\cdots]]-\widetilde{L_k},\]
where:
\[L_k(e)=(\nabla_{Z_1}(\nabla_{Z_2}(\cdots\nabla_{Z_k}(R(X,Y))\cdots)))e.\]
The conclusion follows by observing that $L_k(e)$ can be written in the form:
\[L_k(e)=(\nabla^kR)(Z_1,\ldots,Z_k,X,Y)e+\sum_{i=0}^{k-1}L_{ki},\]
where $L_{ki}$ is a term linear in $(\nabla^iR)({\cdots})e$.
\end{proof}

\end{subsection}

\begin{subsection}{Connections arising from metric tensors}

Let $\pi:E\to M$ be a vector bundle and let $E^*\otimes E^*$ denote the vector bundle over $M$ whose fiber
at $m\in M$ is the space of bilinear forms on $E_m$. If $\nabla$ is a connection on $E$ then we can define
a induced connection $\bil\nabla$ on $E^*\otimes E^*$ by setting:
\[(\bil\nabla_Xg)(\xi,\eta)=X\big(g(\xi,\eta)\big)-g(\nabla_X\xi,\eta)-g(\xi,\nabla_X\eta),\]
where $X$ is a smooth vector field on $M$ and $\xi$, $\eta$ are smooth sections of $E$.
A straightforward computation shows that the curvature tensor $\bil R$ of $\bil\nabla$ is given by:
\begin{equation}\label{eq:Rbil}
\big(\bil R(X,Y)g\big)(\xi,\eta)=-g\big(R(X,Y)\xi,\eta\big)-g\big(\xi,R(X,Y)\eta\big),
\end{equation}
for any smooth vector fields $X$, $Y$ on $M$, any smooth sections $\xi$, $\eta$ of $E$ and any smooth
section $g$ of $E^*\otimes E^*$. If $\gamma:I\to M$ is a smooth curve defined on an interval $I$ around $0$
and if $g_0$ is a bilinear form on $E_{\gamma(0)}$ then the parallel transport $I\ni t\mapsto g_t$ of $g_0$
along $\gamma$ relatively to the connection $\bil\nabla$ is given by:
\[g_t(\xi,\eta)=g_0(P_t^{-1}\xi,P_t^{-1}\eta),\quad\xi,\eta\in E_{\gamma(t)},\]
where $P_t:E_{\gamma(0)}\to E_{\gamma(t)}$ denotes the parallel transport along $\gamma$.

Given a smooth manifold $M$ then a {\em semi-Riemannian metric\/} on $M$ is a smooth section $g$ of the vector bundle
$TM^*\otimes TM^*$ such that $g_m:T_mM\times T_mM\to\R$ is symmetric and nondegenerate; if $g_m$ is positive definite
for all $m\in M$, we call $g$ a {\em Riemannian metric}. The {\em Levi-Civita connection\/} of $g$ is the unique symmetric
connection $\nabla$ on $TM$ such that $\bil\nabla g=0$.

We consider the following problem: {\em given a symmetric
connection $\nabla$ on a smooth manifold $M$, when does there exist a semi-Riemannian metric $g$ on $M$ such that
$\nabla$ is the Levi-Civita connection of $g$}?

Note that if $\nabla$ is the Levi-Civita connection of a semi-Riemannian metric $g$ then for any $m\in M$
and any $v,w\in T_mM$, the linear operator $R_m(v,w):T_mM\to T_mM$ corresponding to the curvature tensor of $\nabla$
is anti-symmetric with respect to $g_m$; moreover, given a smooth curve $\gamma:[a,b]\to M$ with $\gamma(a)=m_0$
and $\gamma(b)=m$ then, denoting by $P:T_{m_0}M\to T_mM$ the parallel transport along $\gamma$, the linear operator:
\[P^{-1}\big[R_m(v,w)\big]P:T_{m_0}M\longrightarrow T_{m_0}M\]
is anti-symmetric with respect to $g_{m_0}$, for all $v,w\in T_mM$. We will show below that this anti-symmetry characterizes
the connections arising from semi-Riemann\-ian metrics.
\begin{prop}\label{thm:vemdemetrica}
Let $M$ be a smooth manifold, $\nabla$ be a symmetric connection on $TM$, $m_0\in M$ and $g_0$ be a nondegenerate
symmetric bilinear form on $T_{m_0}M$. Let $\psi:Z\subset\R\times\Lambda\to M$ be
a $\Lambda$-parametric family of curves on $M$ with a local right inverse $\alpha:V\subset M\to Z$; assume
that $\psi(0,\lambda)=m_0$, for all $\lambda\in M$. For each $(t,\lambda)\in Z$, we denote by
$P_{(t,\lambda)}:T_{m_0}M\to T_{\psi(t,\lambda)}M$ the parallel transport along $t\mapsto\psi(t,\lambda)$.
Assume that for all $(t,\lambda)\in Z$ the linear operator:
\begin{equation}\label{eq:hipoteseR}
P_{(t,\lambda)}^{-1}\big[R_{\psi(t,\lambda)}(v,w)\big]P_{(t,\lambda)}:T_{m_0}M\longrightarrow T_{m_0}M
\end{equation}
is anti-symmetric with respect to $g_0$, for all $v,w\in T_{\psi(t,\lambda)}M$,
where
\[R_{\psi(t,\lambda)}(v,w):T_{\psi(t,\lambda)}M\longrightarrow T_{\psi(t,\lambda)}M\]
denotes the linear operator corresponding to the curvature tensor of $\nabla$. Then $\nabla$ is the Levi-Civita
connection of the semi-Riemannian metric $g$ on $V\subset M$ defined by setting:
\[g_m(\cdot,\cdot)=g_0(P_{\alpha(m)}^{-1}\cdot,P_{\alpha(m)}^{-1}\cdot),\]
for all $m\in V$.
\end{prop}
\begin{proof}
For each $(t,\lambda)\in Z$, let $\tilde\psi(t,\lambda)\in TM^*\otimes TM^*$ be the bilinear form on $T_{\psi(t,\lambda)}M$
defined by:
\[\tilde\psi(t,\lambda)(\cdot,\cdot)=g_0(P_{(t,\lambda)}^{-1}\cdot,P_{(t,\lambda)}^{-1}\cdot).\]
Then $\tilde\psi$ satisfies the hypotheses of Corollary~\ref{thm:corparalelo} with $E=TM^*\otimes TM^*$;
namely, $\tilde\psi(0,\lambda)=g_0$, for all $\lambda\in\Lambda$ and by \eqref{eq:Rbil} and the
anti-symmetry of \eqref{eq:hipoteseR}, we have $\bil R_{\psi(t,\lambda)}(v,w)=0$, for all $v,w\in T_{\phi(t,\lambda)}M$.
Hence $g=\tilde\psi\circ\alpha:V\to TM^*\otimes TM^*$ is a parallel section of $TM^*\otimes TM^*$ and $\nabla$ is
the Levi-Civita connection of $g$.
\end{proof}

\begin{teo}
Let $M$ be a smooth manifold, $\nabla$ be a symmetric connection on $TM$, $m_0\in M$ and $g_0$ be a nondegenerate
symmetric bilinear form on $T_{m_0}M$. Let $\mathcal S$ be a fixed spray on $M$. Assume that:
\begin{itemize}
\item for every piecewise solution $\gamma:[a,b]\to M$ of $\mathcal S$ with $\gamma(a)=m_0$ the linear
operator $P_\gamma^{-1}R_{\gamma(b)}P_\gamma$ on $T_{m_0}M$ is $g_0$-anti-symmetric, where $P_\gamma:T_{m_0}M\to T_{\gamma(b)}M$
denotes parallel transport along $\gamma$;
\item $M$ is (connected and) simply-connected.
\end{itemize}
Then $g_0$ extends to a Riemannian metric on $M$ for which $\nabla$ is the Levi-Civita connection.
\end{teo}
\begin{proof}
It follows from \eqref{eq:Rbil} and Corollary~\ref{thm:corparaleloglobal}.
\end{proof}

\begin{prop}\label{thm:propmetricarealanalitica}
Let $M$ be a (connected and) simply-connected real-analytic manifold and let $\nabla$ be a real-analytic
symmetric connection on $TM$. If there exists a semi-Riemannian metric $g$ on a nonempty open connected subset
of $M$ having $\nabla$ as its Levi-Civita connection then $g$ extends to a globally defined semi-Riemannian metric
on $M$ having $\nabla$ as its Levi-Civita connection.
\end{prop}
\begin{proof}
It follows from Corollary~\ref{thm:correalanalyticparallel}.
\end{proof}

\begin{prop}
Let $M$ be a real-analytic manifold and let $\nabla$ be a real-analytic symmetric connection on $TM$.
Given a point $x_0\in M$ and a nondegenerate symmetric bilinear form $g_0$ on $T_{x_0}M$, if:
\[(\nabla^kR)(v_1,\ldots,v_{k+2}):T_xM\to T_xM\]
is $g_0$-anti-symmetric for all $v_1,\ldots,v_{k+2}\in T_{x_0}M$ and all $k\ge0$ then $g_0$ extends to a semi-Riemannian
metric on an open neighborhood of $x_0$ whose Levi-Civita connection is $\nabla$. Moreover, if $M$ is (connected and)
simply-connected then $g_0$ extends to a global semi-Riemannian metric on $M$ having $\nabla$ as its Levi-Civita connection.
\end{prop}
\begin{proof}
Follows easily from Proposition~\ref{thm:prophighorderhorizontal} and from formula \eqref{eq:Rbil}.
\end{proof}

The above characterizations of Levi--Civita connections have been used in \cite{Cordoba},
where the authors study left-invariant (symmetric) connections in Lie groups.

\end{subsection}



\end{section}

\begin{section}{Affine maps}\label{sec:affine}

Let us now discuss as an application of the ``single leaf Frobenius Theorem'' a classical result in differential geometry.

\subsection{The Cartan--Ambrose--Hicks Theorem}

Consider the following setup. Let $M$, $N$ be smooth manifolds endowed respectively with connections $\nabla^M$ and $\nabla^N$. We denote
by $T^M$, $T^N$ (resp., $R^M$, $R^N$) respectively the torsion tensors (resp., curvature tensors) of $\nabla^M$ and $\nabla^N$.
A smooth map $f:M\to N$ is called {\em affine\/} if for every $x\in M$, $v\in T_xM$ and every smooth vector field $X$ on $M$ we have:
\[\dd f_x(\nabla^M_vX)=\nabla^N_v(\dd f\circ X);\]
in the formula above $\dd f\circ X:M\to TN$ is regarded as a
vector field along $f$ on $N$, so that it makes sense to compute its covariant derivative
$\nabla^N$ along $v\in TM$.

Let $x_0\in M$, $y_0\in N$ be given and let $\sigma_0:T_{x_0}M\to T_{y_0}N$ be a linear map.
Given a geodesic $\gamma:[a,b]\to M$ with $\gamma(a)=x_0$ then the geodesic $\mu:[a,b]\to N$ with $\mu(a)=y_0$ and
$\mu'(a)=\sigma\big(\gamma'(a)\big)$ is called {\em induced\/} on $N$ by the geodesic $\gamma$ and by $\sigma_0$. We observe that
the geodesic $\mu:[a,b]\to N$ is well-defined only if $(b-a)\sigma\big(\gamma'(a)\big)$ is in the domain of the exponential map
of $N$ at the point $y_0$. Let $\sigma:T_{\gamma(b)}M\to T_{\mu(b)}N$ be the linear map given by the composition of parallel transport
along $\gamma$, $\sigma_0$ and parallel transport along $\mu$; we call $\sigma$ the linear map {\em induced\/} by $\gamma$ and $\sigma_0$.

\begin{teo}\label{thm:localCAH}
Let $x_0\in M$, $y_0\in N$ be given and let $\sigma_0:T_{x_0}M\to T_{y_0}N$ be a linear map. Let $U$ be an open subset of $T_{x_0}M$ which is star-shaped at the origin and which is carried diffeomorphically onto an open subset $V$ of $M$ by the exponential map of $M$ at $x_0$.
Assume that $\sigma(U)$ is contained in the domain of the exponential map of $N$ at $y_0$. For each $x\in V$, let $\gamma_x:[0,1]\to M$ be
the unique geodesic such that $\gamma_x'(0)\in U$ and $\gamma_x(1)=x$; let $\mu_x:[0,1]\to N$ and $\sigma_x:T_xM\to T_{\mu_x(1)}N$ be respectively
the geodesic and the linear map induced by $\gamma_x$ and $\sigma_0$. Assume that for all $x\in V$ the linear map $\sigma_x$ relates
$T^M$ with $T^N$ and $R^M$ with $R^N$, i.e.:
\[\sigma_x\big(T^M(\cdot,\cdot)\big)=T^N\big(\sigma_x(\cdot),\sigma_x(\cdot)\big),\quad
\sigma_x\big(R^M(\cdot,\cdot)\cdot\big)=R^N\big(\sigma_x(\cdot),\sigma_x(\cdot)\big)\sigma_x(\cdot).\]
Then the smooth map $f:V\to N$ defined by $f(x)=\mu_x(1)$ is affine and
$\dd f(x)=\sigma_x$ for all $x\in V$; in particular, $f(x_0)=y_0$ and $\dd f(x_0)=\sigma_0$.
\end{teo}

\begin{rem}
In the statement of Theorem~\ref{thm:localCAH}, if one assumes that $\sigma_0$ is an isomorphism (resp., injective) then it follows that $f$ is a local diffeomorphism (resp., that $f$ is an immersion). Moreover, if $\nabla^M$ and $\nabla^N$ are the Levi-Civita connections of Riemannian metrics
on $M$ and $N$ respectively then, if one assumes that $\sigma_0$ is an isometry, it follows that $f$ is a local isometry.
\end{rem}

In what follows we assume that $\nabla^N$ is {\em geodesically complete\/}, i.e., for all $y\in N$ the exponential map of $N$ at $y$ is defined
on the whole tangent space $T_yN$.

Let $x_0\in M$, $y_0\in N$ be given and let $\sigma_0:T_{x_0}M\to T_{y_0}N$ be a linear map. Let $\gamma:[a,b]\to M$
be a {\em piecewise geodesic\/} with $\gamma(a)=x_0$, i.e., there exists a partition $a=t_0<t_1<\cdots<t_k=b$ of $[a,b]$ such that
$\gamma\vert_{[t_i,t_{i+1}]}$ is a geodesic for all $i$. Using the linear map $\sigma_0$ it is possible to define a piecewise geodesic
$\mu:[a,b]\to N$ and a linear map $\sigma:T_{\gamma(b)}M\to T_{\mu(b)}N$ {\em induced\/} by $\gamma$ in the following way: we first define inductively a sequence of geodesics $\mu_i:[t_i,t_{i+1}]\to N$ and of linear maps $\sigma_i:T_{\gamma(t_i)}M\to T_{\mu_i(t_i)}N$. Let $\mu_0$ and $\sigma_1$
be respectively the geodesic and the linear map induced by the geodesic $\gamma\vert_{[t_0,t_1]}$ and by $\sigma_0$. Assuming that
$\mu_i$ and $\sigma_{i+1}$ are defined we let $\mu_{i+1}$ and $\sigma_{i+2}$ be respectively the geodesic and the linear map induced by the geodesic
$\gamma\vert_{[t_i,t_{i+1}]}$ and by $\sigma_{i+1}$. Finally, we let $\mu:[a,b]\to N$ be the piecewise geodesic such that $\mu\vert_{[t_i,t_{i+1}]}=\mu_i$ for all $i$ and we let $\sigma=\sigma_k$.

\begin{teo}[Cartan--Ambrose--Hicks]\label{thm:CAH}
Let $M$, $N$ be smooth manifolds endowed respectively with connections $\nabla^M$ and $\nabla^N$; assume that $\nabla^N$ is geodesically complete
and that $M$ is connected and simply-connected. Let $x_0\in M$, $y_0\in N$ be given and let $\sigma_0:T_{x_0}M\to T_{y_0}N$ be a linear map.
For each piecewise geodesic $\gamma:[a,b]\to M$ with $\gamma(a)=x_0$ denote by $\mu_\gamma:[a,b]\to N$ and by $\sigma_\gamma:T_{\gamma(b)}M\to
T_{\mu_\gamma(b)}N$ respectively the piecewise geodesic and the linear map induced by the piecewise geodesic $\gamma$ and by $\sigma_0$.
Assume that for every piecewise geodesic $\gamma$ the linear map $\sigma_\gamma$ relates $T^M$ with $T^N$ and $R^M$ with $R^N$. Then there exists
a smooth affine map $f:M\to N$ such that for every piecewise geodesic $\gamma:[a,b]\to M$ we have $f\circ\gamma=\mu_\gamma$ and
$\dd f\big(\gamma(b)\big)=\sigma_\gamma$; in particular, $f(x_0)=y_0$ and $\dd f(x_0)=\sigma_0$.
\end{teo}

\begin{rem}
In the statement of the Cartan--Ambrose--Hicks Theorem, if one assumes in addition that $\sigma_0$ is an isomorphism,
and that $\nabla^M$ is geodesically complete then it follows that the affine map $f:M\to N$ is a covering map.
\end{rem}

\begin{cor}
Let $(M,g^M)$, $(N,g^N)$ be Riemannian manifolds with $(N,g^N)$ complete and $M$ connected and simply-connected.
Let $x_0\in M$, $y_0\in N$ be given and let $\sigma_0:T_{x_0}M\to T_{y_0}N$ be a linear isometry onto a subspace of $T_{y_0}N$.
For each piecewise geodesic $\gamma:[a,b]\to M$ with $\gamma(a)=x_0$ denote by $\mu_\gamma:[a,b]\to N$ and by $\sigma_\gamma:T_{\gamma(b)}M\to
T_{\mu_\gamma(b)}N$ respectively the piecewise geodesic and the linear map induced by the piecewise geodesic $\gamma$ and by $\sigma_0$.
Assume that for every piecewise geodesic $\gamma$ the linear map $\sigma_\gamma$ relates $R^M$ with $R^N$. Then there exists
a totally geodesic isometric immersion $f:M\to N$ with $f(x_0)=y_0$ and $\dd f(x_0)=\sigma_0$.
\end{cor}
\begin{proof}
It follows immediately from Theorem~\ref{thm:CAH}; observe that the condition that $f$ is totally geodesic follows from the
fact that $f$ is affine.
\end{proof}

We now show how the proof of Theorems~\ref{thm:localCAH} and \ref{thm:CAH} can be obtained as an application of the local and the global version
of the ``single leaf Frobenius Theorem'' (Theorems~\ref{thm:Frobenius} and \ref{thm:globalFrobenius}).

Consider the vector bundle $E=\Lin(TM,TN)$ over $M\times N$ whose fiber at a point $(x,y)\in M\times N$ is the space of linear maps
$\Lin(T_xM,T_yN)$. Notice that $E$ coincides with the tensor bundle $\pi_1^*(TM^*)\otimes\pi_2^*(TN)$, where $\pi_1$ and $\pi_2$ denote
the projections of the product $M\times N$. The connections $\nabla^M$ and $\nabla^N$ naturally induce a connection $\nabla$ on $E$ given by:
\begin{equation}\label{eq:indconnLinTmTn}(\nabla_{(v,w)}\sigma)(X)=\nabla^N_{(v,w)}\big(\sigma(X)\big)-\sigma(\nabla^M_vX),\end{equation}
where $v\in TM$, $w\in TN$, $X$ is a smooth vector field on $M$ and $\sigma:M\times N\to E$ is a smooth section of $E$.
In the formula above,
$\sigma(X):M\times N\to TN$ is regarded as vector field along the projection $\pi_2:M\times N\to N$ on $N$.

Given a smooth map $f:U\to N$ defined on an open subset $U$ of $M$ then the differential $\dd f:U\to E$ can be regarded as section of $E$ along
the map $U\ni x\mapsto\big(x,f(x)\big)\in M\times N$, so that it makes sense to consider the covariant derivative of $\dd f$ with respect to the
connection $\nabla$.
\begin{lem}\label{thm:affineparallel}
A smooth map $f:U\to N$ defined on an open subset of $M$ is affine if and only if $\dd f$ is parallel with respect to $\nabla$.
\end{lem}
\begin{proof}
Given $v\in TM$ and a smooth vector field $X$ on $U$ we compute:
\[\big(\nabla_v(\dd f)\big)(X)=\nabla^N_v\big(\dd f(X)\big)-\dd f(\nabla^M_vX).\]
The conclusion follows.
\end{proof}

\begin{lem}\label{thm:lemasigmaparallel}
Let $\lambda:t\mapsto\big(\gamma(t),\mu(t),\sigma(t)\big)$ be a smooth curve on $E$, i.e., $\gamma$ is a curve on $M$, $\mu$ is a curve on $N$
and $\sigma(t)$ is a linear map from $T_{\gamma(t)}M$ to $T_{\mu(t)}N$ for all $t$. Then $\lambda$ is parallel with respect to $\nabla$ (or,
equivalently, $\lambda$ is tangent to the horizontal distribution corresponding to $\nabla$) if and only if the following condition holds:
for every $\nabla^M$-parallel vector field $t\mapsto v(t)\in TM$ along $\gamma$, the vector field $t\mapsto\sigma(t)v(t)\in TN$ along
$\mu$ is $\nabla^N$-parallel.
\end{lem}
\begin{proof}
Let $t\mapsto v(t)$ be a vector field along $\gamma$. Let us denote by $\Ddt{}$, $\Ddt M\,$ and $\Ddt N\,$
respectively the covariant derivatives with respect to the parameter $t$ corresponding to the connections $\nabla$, $\nabla^M$ and $\nabla^N$.
The conclusion follows easily from the following formula:
\[\Ddt N\;[\sigma(t)v(t)]=\Big(\Ddt{}\,\sigma(t)\Big)v(t)+\sigma(t)\Ddt M\;\,v(t),\]
observing that $\lambda$ is $\nabla$-parallel if and only if $\Ddt{}\sigma(t)=0$.
\end{proof}

The geometric interpretation of Lemma~\ref{thm:lemasigmaparallel} is given by the following:
\begin{cor}\label{thm:corgeometricparallel}
Let $\lambda$ be as in the statement of Lemma~\ref{thm:lemasigmaparallel} and let $t_0$ in the domain of $\lambda$ be fixed. Then $\lambda$
is parallel with respect to $\nabla$ if and only if the following condition holds: for all $t$, the linear map $\sigma(t):T_{\gamma(t)}M\to T_{\mu(t)}N$
is given by the composition of $\nabla^M$-parallel transport along $\gamma$, $\sigma(t_0)$ and $\nabla^N$-parallel transport along $\mu$.\qed
\end{cor}

We now explain in which form the ``single leaf Frobenius Theorem'' (Theorem~\ref{thm:Frobenius}) is going to be applied. We consider the smooth
submersion $\pi:E\to M$ given by the composition of the canonical projection $E\to M\times N$ with the first projection $\pi_1:M\times N\to M$.
Given $x\in M$, $y\in N$, $\sigma\in\Lin(T_xM,T_yN)$ then the tangent space $T_\sigma E$ is identified with the direct sum of
$T_xM\oplus T_yN$ (the horizontal space corresponding to the connection $\nabla$) and $\Lin(T_xM,T_yN)$ (the tangent space to the fiber).
We will now define a distribution $\mathcal D$ on the manifold $E$ that is horizontal with respect to the submersion $\pi:E\to M$. We set:
\begin{equation}\label{eq:defDCAH}
\mathcal D_\sigma=\Gr(\sigma)\oplus\{0\}\subset T_xM\oplus T_yN\oplus\Lin(T_xM,T_yN)\cong T_\sigma E,
\end{equation}
where $\Gr(\sigma)\subset T_xM\oplus T_yN$ denotes the graph of the linear map $\sigma$.

\begin{lem}\label{thm:affinehorizontal}
Let $s:U\to E$ be a smooth section of $E$ defined on an open subset $U$ of $M$; we write $s(x)=\big(f(x),\sigma(x)\big)$,
where $f:U\to N$ is a smooth map and $\sigma(x)\in\Lin(T_xM,T_{f(x)}N)$, for all $x\in U$. Then $s$ is $\mathcal D$-horizontal if and only
if $\sigma(x)=\dd f(x)$ for all $x\in U$ and $f$ is affine.
\end{lem}
\begin{proof}
Given $x\in U$, $v\in T_xM$ then the component of $\dd s_x(v)$ in $T_xM\oplus T_{f(x)}N$ is equal to $\big(v,\dd f_x(v)\big)$.
Thus, $s$ is $\mathcal D$-horizontal if and only if $\sigma$ is $\nabla$-parallel and $\sigma(x)=\dd f(x)$, for all $x\in U$. The conclusion
follows from Lemma~\ref{thm:affineparallel}.
\end{proof}

\begin{lem}\label{thm:LemaHorizontalLift}
Let $\lambda$ be as in the statement of Lemma~\ref{thm:lemasigmaparallel} and let $t_0$ in the domain of $\lambda$ be fixed.
Assume that $\gamma$ is a geodesic on $M$. Then $\lambda$
is $\mathcal D$-horizontal if and only if the following conditions hold:
\begin{itemize}
\item $\mu$ is a geodesic on $N$;
\item $\mu'(t_0)=\sigma(t_0)\gamma'(t_0)$;
\item for all $t$, the linear map $\sigma(t):T_{\gamma(t)}M\to T_{\mu(t)}N$
is given by the composition of $\nabla^M$-parallel transport along $\gamma$, $\sigma(t_0)$ and $\nabla^N$-parallel transport along $\mu$.
\end{itemize}
\end{lem}
\begin{proof}
Clearly $\lambda$ is $\mathcal D$-horizontal if and only if $\lambda$ is parallel with respect to $\nabla$ and $\mu'(t)=\sigma(t)\gamma'(t)$,
for all $t$. The conclusion follows from Lemma~\ref{thm:lemasigmaparallel} and Corollary~\ref{thm:corgeometricparallel}.
\end{proof}

\begin{cor}
Let $x_0\in M$, $y_0\in N$ be fixed and let $\sigma_0:T_{x_0}M\to T_{y_0}N$ be a linear map. Let $\gamma:[a,b]\to M$ be a piecewise geodesic
with $\gamma(a)=x_0$. Then $\lambda:[a,b]\ni t\mapsto\big(\gamma(t),\mu(t),\sigma(t)\big)\in E$ is the horizontal lift of $\gamma$
with $\lambda(a)=(x_0,y_0,\sigma_0)$ if and only if $\mu:[a,b]\to N$ is the piecewise geodesic induced by $\gamma$ and $\sigma_0$ and $\sigma(t)$
is the linear map induced by $\gamma\vert_{[a,t]}$ and $\sigma_0$, for all $t$.\qed
\end{cor}

\begin{lem}\label{thm:curvatureE}
The curvature tensor $R^E$ of the connection $\nabla$ of $E$ is given by:
\[R^E_{(x,y)}\big((v_1,w_1),(v_2,w_2)\big)\sigma=R^N_y(w_1,w_2)\circ\sigma-\sigma\circ R^M_x(v_1,v_2),\]
for all $(x,y)\in M\times N$, $v_1,v_2\in T_xM$, $w_1,w_2\in T_yN$, $\sigma\in\Lin(T_xM,T_yN)$.\qed
\end{lem}

\begin{lem}\label{thm:lemanablaPQ}
Let $P$, $Q$ be smooth manifolds, $\nabla$ a connection on $Q$ and $h:P\to Q$ be a smooth map. Given smooth vector
fields $X$, $Y$ in $P$ then:
\[\nabla_X\big(\dd h(Y)\big)-\nabla_Y\big(\dd h(X)\big)-\dd h\big([X,Y]\big)=T\big(\dd h(X),\dd h(Y)\big),\]
where $T$ denotes the torsion of $\nabla$.
\end{lem}
\begin{proof}
It is a standard computation in calculus with connections (see Proposition~\ref{thm:RTpullback}).
\end{proof}

We will now compute the Levi form of the distribution $\mathcal D$. Given $x\in M$, $y\in N$, $\sigma\in\Lin(T_xM,T_yN)$,
the Levi form of $\mathcal D$ at the point $\sigma\in E$ is a bilinear map
$\Levi^{\mathcal D}_\sigma:\mathcal D_\sigma\times\mathcal D_\sigma\to T_\sigma E/\mathcal D_\sigma$. We identify
the space $\mathcal D_\sigma$ with $T_xM$ by the isomorphism:
\[T_xM\ni v\longmapsto\big(v,\sigma(v),0\big)\in\mathcal D_\sigma\subset T_xM\oplus T_yN\oplus\Lin(T_xM,T_yN)\cong T_\sigma E.\]
Moreover, the surjective linear map:
\begin{equation}\label{eq:identificaquociente}
T_\sigma E\ni(v,w,\tau)\longmapsto(w-\sigma(v),\tau)\in T_yN\oplus\Lin(T_xM,T_yN)
\end{equation}
has kernel $D_\sigma$ and thus induces an isomorphism from the space $T_\sigma E/\mathcal D_\sigma$ onto $T_yN\oplus\Lin(T_xM,T_yN)$.
Hence, the Levi form of $\mathcal D$ at $\sigma$ will be identified with a bilinear map:
\[\Levi^{\mathcal D}_\sigma:T_xM\times T_xM\longrightarrow T_yN\oplus\Lin(T_xM,T_yN).\]
We now compute $\Levi^{\mathcal D}_\sigma$.
\begin{lem}\label{thm:LemaLeviD}
Given $x\in M$, $y\in N$, $\sigma\in\Lin(T_xM,T_yN)$, the Levi form of $\mathcal D$ at the point $\sigma\in E$ is given by:
\begin{multline*}
\Levi^{\mathcal D}_\sigma(v_1,v_2)
=\Big(\sigma\big(T^M(v_1,v_2)\big)-T^N\big(\sigma(v_1),\sigma(v_2)\big),\\
\sigma\circ R^M_x(v_1,v_2)-R^N_y\big(\sigma(v_1),\sigma(v_2)\big)\circ\sigma\Big),
\end{multline*}
for all $v_1,v_2\in T_xM$.
\end{lem}
\begin{proof}
Given a smooth vector field $X$ on $M$, we define a smooth vector field $\widetilde X$ on $E$ by setting:
\begin{equation}\label{eq:defwidetildeX}
\widetilde X(x,y,\sigma)=\big(X(x),\sigma(X(x)),0\big)\in T_xM\oplus T_yN\oplus\Lin(T_xM,T_yN)\cong T_\sigma E,
\end{equation}
for all $x\in M$, $y\in N$, $\sigma\in\Lin(T_xM,T_yN)$. Observe that $\widetilde X$ is $\mathcal D$-horizontal.

Let $x\in M$, $y\in N$, $\sigma\in\Lin(T_xM,T_yN)$, $v_1,v_2\in T_xM$ be fixed. Choose smooth vector fields $X_1$, $X_2$
on $M$ with $X_1(x)=v_1$, $X_2(x)=v_2$.
In order to compute the Levi form of $\mathcal D$ at the point $\sigma$ it suffices
to compute the Lie bracket $[\widetilde X_1,\widetilde X_2]$ at the point $\sigma$. The vector $[\widetilde X_1,\widetilde X_2]_\sigma$
is identified with an element of $T_xM\oplus T_yN\oplus\Lin(T_xM,T_yN)$. The component in $\Lin(T_xM,T_yN)$ of
such vector can be computed using Lemma~\ref{thm:Levihorizontal}, since $\widetilde X_1$ and $\widetilde X_2$ are
both horizontal with respect to the connection $\nabla$ of $E$; thus, the component of $[\widetilde X_1,\widetilde X_2]_\sigma$
in $\Lin(T_xM,T_yN)$ is equal to $-R^E\big((v_1,\sigma(v_1)),(v_2,\sigma(v_2))\big)\sigma$. Let us now compute
the component of $[\widetilde X_1,\widetilde X_2]_\sigma$ in $T_xM\oplus T_yN$; this is just
$\dd\pi_\sigma\big([\widetilde X_1,\widetilde X_2]_\sigma\big)$. Consider the connection $\nabla^{M\times N}$
on $M\times N$ induced from $\nabla^M$ and $\nabla^N$; its torsion $T^{M\times N}$ is given by:
\[T^{M\times N}\big((v_1,w_1),(v_2,w_2)\big)=\big(T^M(v_1,v_2),T^N(w_1,w_2)\big).\]
We now compute $\dd\pi_\sigma\big([\widetilde X_1,\widetilde X_2]_\sigma\big)$ using Lemma~\ref{thm:lemanablaPQ}
with $P=E$, $Q=M\times N$ and $h=\pi$. We get:
\begin{multline}\label{eq:contachata1}
\nabla^{M\times N}_{\widetilde X_1}\big(\dd\pi(\widetilde X_2)\big)-\nabla^{M\times N}_{\widetilde X_2}\big(\dd\pi(\widetilde X_1)\big)
-\dd\pi\big([\widetilde X_1,\widetilde X_2]\big)\\
=\big(T^M(X_1,X_2),T^N(\sigma(X_1),\sigma(X_2))\big).
\end{multline}
We compute $\nabla^{M\times N}_{\widetilde X_1}\big(\dd\pi(\widetilde X_2)\big)$ as follows:
\[\nabla^{M\times N}_{\widetilde X_1}\big(\dd\pi(\widetilde X_2)\big)=\Ddt{M\times N}\qquad\dd\pi\big(\widetilde X_2(\lambda(t))\big),\]
where $\lambda:\left]-\varepsilon,\varepsilon\right[\to E$ is an integral curve of $\widetilde X_1$ with $\lambda(0)=\sigma$.
Thus $\lambda(t)=\big(x(t),y(t),\sigma(t)\big)$, where $t\mapsto x(t)\in M$ is an integral curve of $X_1$,
$y'(t)=\sigma(t)x'(t)$ and $t\mapsto\sigma(t)$ is $\nabla$-parallel. Hence:
\begin{multline*}
\Ddt{M\times N}\qquad\dd\pi\big(\widetilde X_2(\lambda(t))\big)=\Ddt{M\times N}\quad\ \ \big(X_2(x(t)),\sigma(t)X_2(x(t))\big)\\
=\big(\Ddt M\;X_2(x(t)),\Ddt N\;[\sigma(t)X_2(x(t))]\big)\\
\stackrel{\text{$\sigma$ parallel}}=
\big(\Ddt M\;X_2(x(t)),\sigma(t)\Ddt M\;X_2(x(t))\big).
\end{multline*}
Evaluating at $t=0$ we obtain:
\begin{equation}\label{eq:contachata2}
\nabla^{M\times N}_{\widetilde X_1}\big(\dd\pi(\widetilde X_2)\big)
=\Ddt{M\times N}\qquad\dd\pi\big(\widetilde X_2(\lambda(t))\big)\Big\vert_{t=0}=\big(\nabla^M_{X_1}X_2,
\sigma(\nabla^M_{X_1}X_2)\big),
\end{equation}
where the righthand side of \eqref{eq:contachata2} is evaluated at the point $x$.
Similarly:
\begin{equation}\label{eq:contachata3}
\nabla^{M\times N}_{\widetilde X_2}\big(\dd\pi(\widetilde X_1)\big)=
\big(\nabla^M_{X_2}X_1,\sigma(\nabla^M_{X_2}X_1)\big).
\end{equation}
Using \eqref{eq:contachata1}, \eqref{eq:contachata2} and \eqref{eq:contachata3} we get:
\begin{multline*}
\dd\pi_\sigma\big([\widetilde X_1,\widetilde X_2]_\sigma\big)\\
=\Big([X_1,X_2]_x,\sigma([X_1,X_2]_x)+\sigma\big(T^M(X_1,X_2)\big)
-T^N\big(\sigma(X_1),\sigma(X_2)\big)\Big).
\end{multline*}
Hence, recalling Lemma~\ref{thm:curvatureE}:
\begin{multline}\label{thm:bracktildes}
\!\!\!\![\widetilde X_1,\widetilde X_2]_\sigma=\Big([X_1,X_2]_x,\sigma([X_1,X_2]_x)+\sigma\big(T^M(X_1,X_2)\big)
-T^N\big(\sigma(X_1),\sigma(X_2)\big),\\
\sigma\circ R^M_x(X_1,X_2)-R^N_y(\sigma(X_1),\sigma(X_2))\circ\sigma\Big).
\end{multline}
The conclusion follows recalling formula \eqref{eq:identificaquociente} that gives the identification
between $T_\sigma E/\mathcal D_\sigma$ and $T_yN\oplus\Lin(T_xM,T_yN)$.
\end{proof}

\begin{cor}\label{thm:corLevizero}
Given $x\in M$, $y\in N$, $\sigma\in\Lin(T_xM,T_yN)$, then the Levi form of $\mathcal D$ at the point $\sigma\in E$ vanishes
if and only if the linear map $\sigma:T_xM\to T_yN$ relates $T^M$ with $T^N$ and $R^M$ with $R^N$.
\end{cor}
\begin{proof}
It follows from Lemma~\ref{thm:LemaLeviD}.
\end{proof}

\begin{proof}[Proof of Theorems~\ref{thm:localCAH} and \ref{thm:CAH}]
It follows from Theorems~\ref{thm:Frobenius} and \ref{thm:globalFrobenius}, keeping in mind
Examples~\ref{exa:psigeodesic} and \ref{exa:geodesicspray}, Lemmas~\ref{thm:affinehorizontal} and \ref{thm:LemaHorizontalLift}, and
Corollary~\ref{thm:corLevizero}.
\end{proof}

\subsection{Higher order Cartan--Ambrose--Hicks theorem}

Given a tensor field $\tau$ on a manifold endowed with a connection $\nabla$, we denote by $\nabla^{(r)}\tau$
its $r$-th covariant derivative, for $r\ge1$; we set $\nabla^{(0)}\tau=\tau$.

\begin{teo}\label{thm:highCAH}
Let $M$, $N$ be real-analytic manifolds endowed with real-analytic connections $\nabla^M$ and $\nabla^N$,
respectively.
Let $x_0\in M$, $y_0\in N$ be given and let $\sigma_0:T_{x_0}M\to T_{y_0}N$ be a linear map. If for all $r\ge0$
the linear map $\sigma_0$ relates $\nabla^{(r)}T^M_{x_0}$ with $\nabla^{(r)}T^N_{y_0}$ and
$\nabla^{(r)}R^M_{x_0}$ with $\nabla^{(r)}R^N_{y_0}$ then there exists a real-analytic affine map
$f:U\to N$ defined on an open neighborhood $U$ of $x_0$ in $M$ satisfying $f(x_0)=y_0$ and $\dd f(x_0)=\sigma_0$.
\end{teo}
\begin{proof}
We will apply Theorem~\ref{thm:higherorderFrobenius} to the distribution $\mathcal D$ on $E$ defined in \eqref{eq:defDCAH}.
As before, for $x\in M$, $y\in N$, $\sigma\in\Lin(T_xM,T_yN)$, we use the identification:
\[T_\sigma E\cong T_xM\oplus T_yN\oplus\Lin(T_xM,T_yN).\]
Given a smooth vector field $X$ on $M$, we define a $\mathcal D$-horizontal vector field $\widetilde X$ on $E$
as in \eqref{eq:defwidetildeX}. Recall that for $X_1,X_2\in\Gamma(TM)$, the bracket $[\widetilde X_1,\widetilde X_2]$
was computed in the proof of Lemma~\ref{thm:LemaLeviD} (see \eqref{thm:bracktildes}).
By Remark~\ref{thm:highorderref}, the thesis will follow once we show that the iterated brackets:
\begin{equation}\label{eq:iteratedtildeX}
[\widetilde X_{r+1},\ldots,\widetilde X_1]\stackrel{\text{def}}=
[\widetilde X_{r+1},[\widetilde X_r,\ldots,[\widetilde X_2,\widetilde X_1]\cdots]],
\end{equation}
evaluated at $\sigma_0\in E$ are in $\mathcal D_{\sigma_0}$, for all $X_1,\ldots,X_r\in\Gamma(TM)$ and all $r\ge1$.
For $r\ge0$, $x\in M$, $y\in N$, $\sigma\in\Lin(T_xM,T_yN)$ we set:
\begin{gather*}
\begin{aligned}
\mathfrak T_\sigma^{(r)}(X_1,\ldots,X_{r+2})&=\sigma\big(\nabla^{(r)}T^M(X_1,\ldots,X_{r+2})\big)\\
&-\nabla^{(r)}T^N\big(\sigma(X_1),\ldots,\sigma(X_{r+2})\big)\in T_yN,
\end{aligned}\\[5pt]
\begin{aligned}
\mathfrak R_\sigma^{(r)}(X_1,\ldots,X_{r+2})&=\sigma\circ\nabla^{(r)}R^M(X_1,\ldots,X_{r+2})\\
&-\nabla^{(r)}R^N\big(\sigma(X_1),\ldots,\sigma(X_{r+2})\big)\circ\sigma\in\Lin(T_xM,T_yN),
\end{aligned}
\end{gather*}
for all $X_1,\ldots,X_{r+2}\in T_xM$. The hypotheses of the theorem say that $\mathfrak T_{\sigma_0}^{(r)}$
and $\mathfrak R_{\sigma_0}^{(r)}$ vanish for all $r\ge0$. Observe that $\mathfrak T^{(r)}$ and
$\mathfrak R^{(r)}$ are, respectively, sections along the map $\pi:E\to M\times N$ of the vector bundles
over $M\times N$ given by:
\[\Big(\bigotimes_{r+2}(\pi_1^*TM)^*\Big)\otimes\pi_2^*TN\quad\text{and}\quad
\Big(\bigotimes_{r+2}(\pi_1^*TM)^*\Big)\otimes(\pi_1^*TM)^*\otimes\pi_2^*TN,\]
where $\pi_1$ and $\pi_2$ denote the projections of the product $M\times N$.

Our plan is to show that the iterated bracket
\eqref{eq:iteratedtildeX} can be written in the form:
\begin{multline}\label{thm:superformula}
[\widetilde X_{r+1},\ldots,\widetilde X_1]=\big(0,\mathfrak T^{(r-1)}(X_{r+1},\ldots,X_1),
\mathfrak R^{(r-1)}(X_{r+1},\ldots,X_1)\big)\\
+\big(0,\mathcal L^{(r)}(\mathfrak T^{(0)},\mathfrak R^{(0)},\ldots,\mathfrak T^{(r-2)},\mathfrak R^{(r-2)})\big)
+\text{terms in $\Gamma^{(r-1)}(\mathcal D)$},
\end{multline}
for all $r\ge1$, where $\mathcal L^{(r)}$ is a section along the map $\pi:E\to M\times N$ of the vector bundle
over $M\times N$ given by:
\begin{multline*}
\Lin\bigg(\bigoplus_{i=0}^{r-2}\Big[\Big(\bigotimes_{i+2}(\pi_1^*TM)^*\Big)\otimes\pi_2^*TN
\oplus
\Big(\bigotimes_{i+2}(\pi_1^*TM)^*\Big)\otimes(\pi_1^*TM)^*\otimes\pi_2^*TN\Big],\\
\pi_2^*TN\oplus(\pi_1^*TM)^*\otimes\pi_2^*TN\bigg).
\end{multline*}
Once formula \eqref{thm:superformula} is proven, the conclusion follows easily by induction on $r$.
We will now conclude the proof by showing formula \eqref{thm:superformula} by induction on $r$.
For $r=1$, we have (recall \eqref{thm:bracktildes}):
\begin{multline*}
[\widetilde X_2,\widetilde X_1]=\Big(0,\sigma\big(T^M(X_2,X_1)\big)
-T^N\big(\sigma(X_2),\sigma(X_1)\big),\\
\sigma\circ R^M(X_2,X_1)-R^N(\sigma(X_2),\sigma(X_1))\circ\sigma\Big)
+\text{terms in $\Gamma^0(\mathcal D)$}\\
=\big(0,\mathfrak T^{(0)}(X_2,X_1),\mathfrak R^{(0)}(X_2,X_1)\big)+\text{terms in $\Gamma^0(\mathcal D)$},
\end{multline*}
proving the base of the induction. The induction step can be proven by applying $\ad_{\widetilde X_{r+2}}$
to both sides of \eqref{thm:superformula}, keeping in mind Lemma~\ref{thm:LemaAB} below and the following formulas:
\begin{gather*}
(\nabla_{\mathrm{hor}}\mathfrak T^{(i)})_\sigma\big(Z,\sigma(Z)\big)=\mathfrak T^{(i+1)}(Z,\cdots),\\
(\nabla_{\mathrm{hor}}\mathfrak R^{(i)})_\sigma\big(Z,\sigma(Z)\big)=\mathfrak R^{(i+1)}(Z,\cdots),
\end{gather*}
where $Z$ is a vector field on $M$. This concludes the proof.
\end{proof}

\begin{lem}\label{thm:LemaAB}
Let $A$ be a section along the map $\pi_2\circ\pi:E\to N$ of the tangent bundle of $N$ and let
$B$ be a section along the map $\pi:E\to M\times N$ of the vector bundle $E$, so that
$(0,A,B)$ is a vector field on the manifold $E$. Let $Z$ be a vector field on $M$.
Then:
\begin{multline*}
[\widetilde Z,(0,A,B)]_\sigma=\Big(0,\nabla^N_{\mathrm{hor}}A\big(Z,\sigma(Z)\big)-B(Z)-T^N\big(\sigma(Z),A\big),\\
\nabla_{\mathrm{hor}}B\big(Z,\sigma(Z)\big)-R^N\big(\sigma(Z),A\big)\circ\sigma\Big),\quad\sigma\in E,
\end{multline*}
where $\nabla^N_{\mathrm{hor}}A$ (resp., $\nabla_{\mathrm{hor}}B$)
denotes the restriction of $\nabla^NA$ (resp., of $\nabla B$)
to the horizontal subbundle of $TE$ determined by the connection of $E$.
\end{lem}
\begin{proof}
We compute the horizontal component of $[\widetilde Z,(0,A,B)]$ using Lemma~\ref{thm:lemanablaPQ}
with $P=E$, $Q=M\times N$ and $h=\pi$. We have:
\begin{multline*}
\dd\pi_\sigma[(Z,\sigma(Z),0),(0,A,B)]=\nabla^{M\times N}_{(Z,\sigma(Z),0)}(0,A)
-\nabla^{M\times N}_{(0,A,B)}\big(Z,\sigma(Z)\big)\\
-T^{M\times N}\big((Z,\sigma(Z)),(0,A)\big).
\end{multline*}
Clearly:
\[T^{M\times N}\big((Z,\sigma(Z)),(0,A)\big)=\big(T^M(Z,0),T^N(\sigma(Z),A)\big)=\big(0,T^N(\sigma(Z),A)\big)\]
and:
\[\nabla^{M\times N}_{(Z,\sigma(Z),0)}(0,A)=\big(0,\nabla^N_{\mathrm{hor}}A(Z,\sigma(Z))\big).\]
Let $t\mapsto\big(x(t),y(t),\sigma(t)\big)$ be an integral curve of $(0,A,B)$, i.e., $t\mapsto x(t)$ is constant,
$y'=A$ and $\Ddt{}\sigma=B$. We compute:
\[\nabla^{M\times N}_{(0,A,B)}\big(Z,\sigma(Z)\big)=\Ddt{M\times N}\quad\ \ \big(Z_{x(t)},\sigma(t)Z_{x(t)}\big)
=\big(0,B(Z)\big).\]
Let us now compute the vertical component of $[\widetilde Z,(0,A,0)]$. Since both $\widetilde Z$ and $(0,A,0)$
are in the horizontal subbundle of $TE$ determined by the connection of $E$, the vertical component
of $[\widetilde Z,(0,A,0)]$ can be directly computed using Lemmas~\ref{thm:Levihorizontal} and \ref{thm:curvatureE}, as follows:
\begin{multline*}
\text{vertical component of $[\widetilde Z,(0,A,0)]_\sigma$}=-R^E\big((Z,\sigma(Z)),(0,A)\big)\sigma\\
=\sigma\circ R^M(Z,0)-R^N\big(\sigma(Z),A\big)\circ\sigma=-R^N\big(\sigma(Z),A\big)\circ\sigma.
\end{multline*}
Finally, we compute the vertical component of $[\widetilde Z,(0,0,B)]$. Let $W$ be a vector field on $M$ and
$\alpha$ be a $1$-form on $N$; we define a map $f_{W,\alpha}:E\to\R$ by setting:
\[f_{W,\alpha}(\sigma)=\alpha\big(\sigma(W)\big).\]
Let $x\in M$, $y\in N$, $\sigma\in\Lin(T_xM,T_yN)$ be fixed and assume that
$\nabla^MW(x)=0$, $\nabla^N\alpha(y)=0$, so that $\dd f_{W,\alpha}(\sigma)$ annihilates the horizontal subspace
of $T_\sigma E$ determined by the connection of $E$. We compute:
\begin{gather*}
(0,0,B)(f_{W,\alpha})=\alpha\big(B(W)\big),\\
\widetilde Z(f_{W,\alpha})=(\nabla_{\sigma(Z)}\alpha)\big(\sigma(W)\big)+\alpha\big(\sigma(\nabla_ZW)\big),\\
\begin{aligned}
\widetilde Z\big((0,0,B)(f_{W,\alpha})\big)_\sigma
&=(\nabla_{\sigma(Z)}\alpha)_x\big(B(W)\big)+
\alpha\big(\nabla_{\mathrm{hor}}B_\sigma(Z,\sigma(Z))(W)\big)\\
&+\alpha\big(B(\nabla_{Z_x}W)\big)
=\alpha\big(\nabla_{\mathrm{hor}}B_\sigma(Z,\sigma(Z))(W)\big),
\end{aligned}\\
\begin{aligned}
(0,0,B)\big(\widetilde Z(f_{W,\alpha})\big)_\sigma&=
(\nabla_{B(Z)}\alpha)_x\big(\sigma(W)\big)+(\nabla_{\sigma(Z)}\alpha)_x\big(B(W)\big)\\
&+\alpha\big(B(\nabla_{Z_x}W)\big)=0,
\end{aligned}
\end{gather*}
so that:
\[[\widetilde Z,(0,0,B)](f_{W,\alpha})=\alpha\big(\nabla_{\mathrm{hor}}B_\sigma(Z,\sigma(Z))(W)\big).\]
Hence, the vertical component of $[\widetilde Z,(0,0,B)]_\sigma$ is equal to
$\nabla_{\mathrm{hor}}B_\sigma\big(Z,\sigma(Z)\big)$. This concludes the proof.
\end{proof}

\begin{prop}\label{thm:globalRACAH}
Let $M$, $N$ be real-analytic manifolds endowed with real-analytic connections $\nabla^M$ and $\nabla^N$,
respectively. Assume that $\nabla^N$ is geodesically complete and that $M$ is (connected and) simply-connected.
Then every affine map $f:U\to N$ defined on a nonempty connected open subset $U$ of $M$ extends to an affine
map from $M$ to $N$. In particular, if in addition $x_0\in M$, $y_0\in N$, $\sigma_0\in\Lin(T_{x_0}M,T_{y_0}N)$ satisfy the
hypotheses of Theorem~\ref{thm:highCAH} then there exists an affine map $f:M\to N$ with $f(x_0)=y_0$ and
$\dd f(x_0)=\sigma_0$.
\end{prop}
\begin{proof}
If $\mathcal D$ is the distribution on $E$ defined in \eqref{eq:defDCAH} then, by Lemma~\ref{thm:affinehorizontal},
$s(x)=\big(f(x),\dd f(x)\big)$ is a $\mathcal D$-horizontal section of $\pi:E\to M$ defined in $U$.
The geodesical completeness of $\nabla^N$ guarantees that hypothesis~(b) of Proposition~\ref{thm:Frobeniusrealanalytic}
is satisfied; hence, such proposition gives a global horizontal section of $\pi$.
\end{proof}

An {\em affine symmetry\/} around a point $x_0\in M$ is an affine map $f:U\to M$ defined in an open neighborhood
$U$ of $x_0$ with $f(x_0)=x_0$ and $\dd f(x_0)=-\Id$.
\begin{cor}\label{thm:curiouscor}
Let $M$ be a real-analytic manifold endowed with a real-analytic connection $\nabla$. Let $x_0\in M$ be fixed.
Then there exists an affine symmetry around $x_0$ if and only if:
\begin{equation}\label{eq:nablaparimpar}
\nabla^{(2r)}T_{x_0}=0,\quad\text{and}\quad\nabla^{(2r+1)}R_{x_0}=0,\quad
\text{for all $r\ge0$}.
\end{equation}
Moreover, if $M$ is (connected and) simply-connected and complete then condition \eqref{eq:nablaparimpar} is equivalent
to the existence of a globally defined affine symmetry $f:M\to M$ around $x_0$.
\end{cor}

\begin{proof}
Apply Theorem~\ref{thm:highCAH} with $M=N$, $y_0=x_0$ and $\sigma_0=-\Id$. For the global result apply Proposition~\ref{thm:globalRACAH}.
\end{proof}

\end{section}

\appendix

\begin{section}{A globalization principle}
\label{sec:globalization}

\begin{defin}\label{thm:fundopenset}
Let $X$, $\widetilde X$ be topological spaces and $\pi:\widetilde X\to X$ be a map.
An open subset $U\subset X$ is called a {\em fundamental open subset\/} of $X$ if
$\pi^{-1}(U)$ equals a disjoint union $\bigcup_{i\in I}U_i$ of open subsets $U_i$ of $\widetilde X$
such that $\pi\vert_{U_i}:U_i\to U$ is a homeomorphism for all $i\in I$. We say that $\pi$
is a {\em covering map\/} if $X$ can be covered by fundamental open subsets.
\end{defin}
Obviously every covering map is a local homeomorphism.

Given a local homeomorphism $\pi:\widetilde X\to X$ then by a {\em local section\/} of $\pi$ we mean
a continuous map $s:U\to\widetilde X$ defined on an open subset of $X$ with $\pi\circ s=\Id_U$.
\begin{lem}\label{thm:impliesfundopen}
Let $X$, $\widetilde X$ be topological spaces and $\pi:\widetilde X\to X$ be a local homeomorphism.
Assume that $\widetilde X$ is Hausdorff. Let $U$ be a connected open subset of $X$ satisfying the following property:
\begin{itemize}
\item[($*$)] for every $x\in U$ and every $\tilde x\in\widetilde X$ with $\pi(\tilde x)=x$
there exists a local section $s:U\to\widetilde X$ of $\pi$ with $s(x)=\tilde x$.
\end{itemize}
Then $U$ is a fundamental open subset of $X$.
\end{lem}
\begin{proof}
Let $\mathcal S$ be the set of all local sections of $\pi$ defined in $U$. We claim that:
\[\pi^{-1}(U)=\bigcup_{s\in\mathcal S}s(U).\]
Indeed, if $s\in\mathcal S$ then obviously $s(U)\subset\pi^{-1}(U)$; moreover, given
$\tilde x\in\pi^{-1}(U)$ then $x=\pi(\tilde x)\in U$ and by property ($*$) there exists
$s\in\mathcal S$ with $s(x)=\tilde x$. Thus $\tilde x\in s(U)$. This proves the claim.
Now observe that $s(U)$ is open in $\widetilde X$ for all
$s\in\mathcal S$; moreover, $\pi\vert_{s(U)}:s(U)\to U$ is a homeomorphism, being the inverse
of $s:U\to s(U)$. To complete the proof, we show that the union $\bigcup_{s\in\mathcal S}s(U)$
is disjoint. Pick $s,s'\in\mathcal S$ with $s(U)\cap s'(U)\ne\emptyset$. Then there exists
$x,y\in U$ with $s(x)=s'(y)$. Observe that:
\[x=\pi\big(s(x)\big)=\pi\big(s'(y)\big)=y,\]
and thus $s(x)=s'(x)$. Since $U$ is connected and $\widetilde X$ is Hausdorff it follows that $s=s'$.
\end{proof}

\begin{cor}
Let $X$, $\widetilde X$ be topological spaces and $\pi:\widetilde X\to X$ be a local homeomorphism.
Assume that $\widetilde X$ is Hausdorff and that $X$ is locally connected. If every point of $X$ has an open
neighborhood satisfying property ($*$) then $\pi$ is a covering map.\qed
\end{cor}

Let $X$ be a topological space. A {\em pre-sheaf\/} on $X$ is a map $\mathfrak P$
that assigns to each open subset $U\subset X$ a set $\mathfrak P(U)$ and to each
pair of open subsets $U,V\subset X$ with $V\subset U$ a map $\mathfrak P_{U,V}:\mathfrak P(U)\to\mathfrak P(V)$
such that the following properties hold:
\begin{itemize}
\item for every open subset $U\subset X$ the map $\mathfrak P_{U,U}$ is the identity
map of the set $\mathfrak P(U)$;
\item given open sets, $U,V,W\subset X$ with
$W\subset V\subset U$ then:
\[\mathfrak P_{V,W}\circ\mathfrak P_{U,V}=\mathfrak P_{U,W}.\]
\end{itemize}
We say that the pre-sheaf $\mathfrak P$ is {\em nontrivial\/} if there exists a nonempty open subset $U$ of $X$
with $\mathfrak P(U)\ne\emptyset$.

\begin{example}
For each open subset $U$ of $X$ let $\mathfrak P(U)$ be the set of all continuous maps $f:U\to\R$.
Given open subsets $U$, $V$ of $X$ with $V\subset U$ we set $\mathfrak P_{U,V}(f)=f\vert_V$, for all
$f\in\mathfrak P(U)$. Then $\mathfrak P$ is a pre-sheaf over $X$.
\end{example}

A {\em sheaf\/} over a topological space $X$ is a pair $(\mathcal S,\pi)$, where
$\mathcal S$ is a topological space and $\pi:\mathcal S\to X$ is a local homeomorphism.

Let $\mathfrak P$ be a pre-sheaf over a topological space $X$. Given a point $x\in X$,
consider the disjoint union of all sets $\mathfrak P(U)$, where $U$ is an open neighborhood
of $x$ in $X$. We define an equivalence relation $\sim$ on such disjoint union as follows;
given $f_1\in\mathfrak P(U_1)$, $f_2\in\mathfrak P(U_2)$, where $U_1$, $U_2$ are open neighborhoods
of $x$ in $X$ then $f_1\sim f_2$ if and only if there exists an open neighborhood
$V$ of $x$ contained in $U_1\cap U_2$ such that $\mathfrak P_{U_1,V}(f_1)=\mathfrak P_{U_2,V}(f_2)$.
If $U$ is an open neighborhood of $x$ in $X$ and $f\in\mathfrak P(U)$ then the equivalence
class of $f$ corresponding to the equivalence relation $\sim$ will be denote by $[f]_x$
and will be called the {\em germ\/} of $f$ at the point $x$. We set:
\[\mathfrak S_x=\big\{[f]_x:\text{$f\in\mathfrak P(U)$, for some open neighborhood $U$ of $x$ in $X$}\big\}.\]
Let $\mathfrak S$ denote the disjoint union of all $\mathfrak S_x$, with $x\in X$.
Let $\pi:\mathfrak S\to X$ denote the map that carries $\mathfrak S_x$ to the point $x$.
Our goal now is to define a topology on $\mathfrak S$. Given an open subset $U\subset X$
and an element $f\in\mathfrak P(U)$ we set:
\[\mathcal V(f)=\big\{[f]_x:x\in U\big\}\subset\mathfrak S.\]
The set:
\[\big\{\mathcal V(f):\text{$f\in\mathfrak P(U)$, $U$ an open subset of $X$}\big\}\]
is a basis for a topology on $\mathfrak S$; moreover, if $\mathfrak S$ is endowed with such topology, the map $\pi:\mathfrak S\to X$
is a local homeomorphism, so that $(\mathfrak S,\pi)$ is a sheaf over $X$. We call $(\mathfrak S,\pi)$ the {\em sheaf of germs\/}
corresponding to the pre-sheaf $\mathfrak P$.
Observe that if $U$ is an open subset of $X$ and $f\in\mathfrak P(U)$ then the map $\hat f:U\to\mathfrak S$ defined by
\[\hat f(x)=[f]_x,\quad x\in U,\]
is a local section of the sheaf of germs defined in $U$.

\begin{defin}\label{thm:localizationproperty}
We say that the pre-sheaf $\mathfrak P$ has the {\em localization property\/} if,
given a family $(U_i)_{i\in I}$ of open subsets of $X$ and setting $U=\bigcup_{i\in I}U_i$
then the map:
\begin{equation}\label{eq:PUUi}
\mathfrak P(U)\ni f\longmapsto\big(\mathfrak P_{U,U_i}(f)\big)_{i\in I}\in\prod_{i\in I}\mathfrak P(U_i)
\end{equation}
is injective and its image consists of all the families $(f_i)_{i\in I}$ in $\prod_{i\in I}\mathfrak P(U_i)$
such that $\mathfrak P_{U_i,U_i\cap U_j}(f_i)=\mathfrak P_{U_j,U_i\cap U_j}(f_j)$, for all
$i,j\in I$.
\end{defin}

\begin{rem}\label{thm:remlocalization}
If the pre-sheaf $\mathfrak P$ has the localization property then for every local section $s:U\to\mathfrak S$
of its sheaf of germs $\mathfrak S$ there exists a unique $f\in\mathfrak P(U)$ such that $\hat f=s$.
\end{rem}

\begin{defin}\label{thm:uniquenessprop}
We say that the pre-sheaf $\mathfrak P$ has the {\em uniqueness property\/}
if for every connected open subset $U\subset X$ and every nonempty open subset $V\subset U$
the map $\mathfrak P_{U,V}$ is injective.
We say that an open subset $U\subset X$ has the {\em extension property with respect to the
pre-sheaf\/ $\mathfrak P$\/} if for every connected nonempty open subset $V$ of $U$ the map $\mathfrak P_{U,V}$
is surjective. We say that the pre-sheaf $\mathfrak P$ has the {\em extension property\/}
if $X$ can be covered by open sets having the extension property with respect to $\mathfrak P$.
\end{defin}

\begin{rem}\label{thm:uniquenessHausdorffextensionstar}
If the pre-sheaf $\mathfrak P$ has the uniqueness property and if $X$ is locally connected
and Hausdorff then the space $\mathfrak S$ is Hausdorff.
If $X$ is locally connected and if $U$ is an open subset of $X$ having the extension property with respect to the pre-sheaf
$\mathfrak P$ then $U$ has the property ($*$)
with respect to the local homeomorphism $\pi:\mathfrak S\to X$. It follows from Lemma~\ref{thm:impliesfundopen} that
if $X$ is Hausdorff and locally connected and if the pre-sheaf\/ $\mathfrak P$ has the uniqueness property and
the extension property then the map $\pi:\mathfrak S\to X$ is a covering map.
\end{rem}

\begin{prop}\label{thm:niceglobalization}
Assume that $X$ is Hausdorff, locally arc-connected, connected, and simply-connected.
If $\mathfrak P$ is a pre-sheaf over $X$ satisfying the localization property, the uniqueness property and
the extension property then the open set $X$ has the extension property for $\mathfrak P$, i.e.,
for every nonempty open connected subset $V$ of $X$ the map $\mathfrak P_{X,V}:\mathfrak P(X)\to\mathfrak P(V)$
is surjective. In particular, if $\mathfrak P$ is nontrivial then the set $\mathfrak P(X)$ is nonempty.
\end{prop}
\begin{proof}
Let $V$ be a nonempty open connected subset of $X$ and let $f\in\mathfrak P(V)$ be fixed. We will show that
$f$ is in the image of $\mathfrak P_{X,V}$.
Let $\pi:\mathfrak S\to X$ denote the sheaf of germs of $\mathfrak P$. By Remark~\ref{thm:uniquenessHausdorffextensionstar}
$\pi$ is a covering map. Choose an arbitrary point $x_0\in V$ and let $\mathfrak S_0$
be the arc-connected component of $[f]_{x_0}$ in $\mathfrak S$. Since $X$ is locally arc-connected, the restriction
of $\pi$ to $\mathfrak S_0$ is again a covering map. By the connectedness and simply-connectedness of $X$,
$\pi\vert_{\mathfrak S_0}:\mathfrak S_0\to X$ is a homeomorphism. The inverse of $\pi\vert_{\mathfrak S_0}$ is therefore
a global section $s:X\to\mathfrak S$ and, by Remark~\ref{thm:remlocalization}, there exists $g\in\mathfrak P(X)$ with
$\hat g=s$. Now $[g]_{x_0}=\hat g(x_0)=s(x_0)=[f]_{x_0}$ and hence, by the uniqueness property,
$\mathfrak P_{X,V}(g)=f$.
\end{proof}

\end{section}

\begin{section}{A crash course on calculus with connections}
\label{app:crashcourse}
Given a smooth vector bundle $\pi:E\to M$ over a smooth manifold $M$, we will denote by $\Gamma(E)$ the space
of all smooth sections $s:M\to E$ of $E$. Observe that $\Gamma(E)$ is a real vector space and it is a module
over the commutative ring $C^\infty(M)$ of all smooth maps $f:M\to\R$. Given an open subset $U$ of $M$, we denote
by $E\vert_U$ the restriction of the vector bundle $E$ to $U$, i.e., $E\vert_U=\pi^{-1}(U)$.
\begin{defin}
A {\em connection\/} on a vector bundle $\pi:E\to M$ is a $\R$-bilinear map:
\[\nabla:\Gamma(TM)\times\Gamma(E)\ni(X,s)\longmapsto\nabla_Xs\in\Gamma(E)\]
that is $C^\infty(M)$-linear in the variable $X$ and satisfies the Leibnitz derivative rule:
\[\nabla_X(fs)=X(f)s+f\nabla_Xs,\]
for all $X\in\Gamma(TM)$, $s\in\Gamma(E)$, $f\in C^\infty(M)$.
\end{defin}

\begin{example}
If $E_0$ is a fixed real finite-dimensional vector space and $E=M\times E_0$ is a trivial vector bundle
over $M$ then a section $s$ of $E$ can be identified with a map $s:M\to E_0$ and a connection on $E$ can be defined by:
\begin{equation}\label{eq:trivialconnection}
\nabla_Xs=\dd s(X),
\end{equation}
for all $X\in\Gamma(TM)$. We call \eqref{eq:trivialconnection} the {\em standard connection\/} of the trivial bundle $E$.
\end{example}

It follows from the $C^\infty(TM)$-linearity of $\nabla$ in the variable $X$ that $\nabla_Xs(x)$ depends
only of the value of $X$ at the point $x\in M$, i.e., if $X(x)=X'(x)$ then
$\nabla_Xs(x)=\nabla_{X'}s(x)$. Given $s\in\Gamma(E)$, $x\in M$ and $v\in T_xM$, we set:
\[\nabla_vs=\nabla_Xs(x),\]
where $X\in\Gamma(TM)$ is an arbitrary vector field with $X(x)=v$. For all $x\in M$ we denote by
$\nabla s(x):T_xM\to E_x$ the linear map given by $v\mapsto\nabla_vs$. Thus, given $s\in\Gamma(E)$, we obtain
a smooth section $\nabla s$ of $TM^*\otimes E$.

It follows from the Leibnitz rule that if $U\subset M$ is an open subset
then the restriction of $\nabla_Xs$ to $U$ depends only of the restriction of $s$ to $U$. Thus,
given an open subset $U$ of $M$, a connection $\nabla$ on $E$ induces a unique connection $\nabla^U$ on $E\vert_U$ such that:
\begin{equation}\label{eq:nablaU}
\nabla^U_v(s\vert_U)=\nabla_vs,
\end{equation}
for all $s\in\Gamma(E)$, $v\in TU$.

\begin{rem}
Given connections $\nabla$ and $\nabla'$ on a vector bundle $\pi:E\to M$ then their difference is a tensor; more
explicitly:
\[\mathfrak t(X,s)=\nabla_Xs-\nabla'_Xs\in\Gamma(E),\quad X\in\Gamma(TM),\ s\in\Gamma(E),\]
is $C^\infty(M)$-bilinear and hence defines a smooth section $\mathfrak t$ of the vector bundle
$TM^*\otimes E^*\otimes E$. Moreover, if $\nabla$ is a connection on $E$ and $\mathfrak t$ is a smooth section
of $TM^*\otimes E^*\otimes E$ then $\nabla+\mathfrak t$ is also a connection on $E$.
If $\mathfrak t$ is a section of $TM^*\otimes E^*\otimes E$ then, given $x\in M$, $v\in T_xM$, we identify
$\mathfrak t(v)$ with a linear operator on the fiber $E_x$.
\end{rem}

Given vector bundles $\pi:E\to M$, $\tilde\pi:\widetilde E\to M$ over the same base manifold $M$ then a
{\em vector bundle morphism\/} is a smooth map $L:E\to\widetilde E$ such that $\tilde\pi\circ L=\pi$ and
such that $L\vert_{E_x}:E_x\to\widetilde E_x$ is a linear map, for all $x\in M$. We will denote the restriction
of $L$ to $E_x$ by $L_x$. If $L:E\to\widetilde E$ is a vector bundle morphism such that $L_x$ is an isomorphism
for all $x\in M$ then we call $L$ a {\em vector bundle isomorphism}. If $A$ is an open subset of $E$
and $L:A\to\widetilde E$ is a smooth map such that $\tilde\pi\circ L=\pi\vert_A$ then we call $L$ a
{\em fiber bundle morphism}. Given $x\in M$, we write $A_x=A\cap E_x$ and $L_x=L\vert_{A_x}:A_x\to\widetilde E_x$.
\begin{defin}
Given vector bundles $\pi:E\to M$, $\tilde\pi:\widetilde E\to M$ over the same base manifold $M$, a vector
bundle morphism $L:E\to\widetilde E$ and connections $\nabla$, $\widetilde\nabla$ on $E$ and $\widetilde E$ respectively
then we say that $\nabla$ and $\widetilde\nabla$ are {\em $L$-related\/} if:
\[\widetilde\nabla_X\big(L(s)\big)=L(\nabla_Xs),\]
for all $X\in\Gamma(TM)$, $s\in\Gamma(E)$.
\end{defin}

In what follows, we will deal with several constructions involving connections on different vector bundles.
In order to avoid heavy notations, we will usually denote all these connections by the symbol $\nabla$; it should
be clear from the context which connection the symbol $\nabla$ refers to. For instance, formula \eqref{eq:nablaU} will
be rewritten in the following simpler form:
\[\nabla_v(s\vert_U)=\nabla_vs.\]

\begin{defin}
Given a connection $\nabla$ on a vector bundle $\pi:E\to M$, then the {\em curvature tensor\/} of $\nabla$ is defined
by:
\begin{equation}\label{eq:defRnabla}
R(X,Y)s=\nabla_X\nabla_Ys-\nabla_Y\nabla_Xs-\nabla_{[X,Y]}s\in\Gamma(E),
\end{equation}
for all $X,Y\in\Gamma(TM)$, $s\in\Gamma(E)$.
\end{defin}
Since the righthand side of \eqref{eq:defRnabla} is $C^\infty(M)$-linear in the variables $X$, $Y$ and $s$, it follows
that $R$ can be identified with a smooth section of the vector bundle $TM^*\otimes TM^*\otimes E^*\otimes E$.
Clearly, $R(X,Y)s$ is anti-symmetric in the variables $X$ and $Y$.

The notion of torsion is usually defined only for connection on tangent bundles. We will present a slight generalization
of this notion.
\begin{defin}
Let $\pi:E\to M$ be a smooth vector bundle and let $\iota:TM\to E$ be a vector bundle morphism.
Given a connection $\nabla$ on $E$ then the {\em $\iota$-torsion\/} of $\nabla$ is defined by:
\begin{equation}\label{eq:Tiota}
T^\iota(X,Y)=\nabla_X\big(\iota(Y)\big)-\nabla_Y\big(\iota(X)\big)-\iota\big([X,Y]\big)\in\Gamma(E),
\end{equation}
for all $X,Y\in\Gamma(TM)$. If $E=TM$ and $\iota$ is the identity map of $TM$, we will write simply $T$ and call it
the {\em torsion\/} of $\nabla$.
\end{defin}
Again, the righthand side of \eqref{eq:Tiota} is $C^\infty(M)$-linear on the variables $X$ and $Y$, so that $T^\iota$
can be identified with a smooth section of the vector bundle $TM^*\otimes TM^*\otimes E$. Clearly, $T^\iota(X,Y)$ is anti-symmetric
in $X$ and $Y$.

In what follows we will study some natural constructions with vector bundles endowed with connections and we will
present some formulas for the computation of torsions and curvatures. We will consider constructions that
act on the basis of the vector bundles and constructions that act on their fibers.

Given smooth manifolds, $M$, $N$, a smooth vector bundle $\pi:E\to M$ over $M$ and a smooth map $f:N\to M$
then we denote by $f^*E$ the {\em pull-back of $E$ by $f$\/} which is a vector bundle over $N$ whose fiber
at a point $x\in N$ is equal to $E_{f(x)}$. Observe that there is a natural identification of smooth
sections of the bundle $f^*E$ with smooth {\em sections of $E$ along $f$}, i.e., smooth maps $s:N\to E$ such that
$\pi\circ s=f$. Notice that every smooth section $s:M\to E$ of $E$ gives rise to a smooth section of $E$ along $f$ given
by $s\circ f:N\to E$; we may thus identify $s\circ f$ with a section of $f^*E$.

\begin{prop}\label{thm:defpullbackconnection}
Given smooth manifolds $M$, $N$, a smooth vector bundle $E$ over $M$ endowed
with a connection $\nabla$ and a smooth map $f:N\to M$
then there exists a unique connection $f^*\nabla$ on the pull-back
bundle $f^*E$ such that:
\begin{equation}\label{eq:defnablapull}
(f^*\nabla)_v(s\circ f)=\nabla_{\dd f(v)}s,
\end{equation}
for all $s\in\Gamma(E)$ and all $v\in TN$.
\end{prop}

The next result follows easily from Proposition~\ref{thm:defpullbackconnection}.
\begin{prop}
Let $P$, $N$, $M$ be smooth manifolds, $E$ be a vector bundle over $M$ endowed with
a connection $\nabla$ and $g:P\to N$, $f:N\to M$ be smooth
maps. Then:
\begin{equation}\label{eq:pullpull}
(f\circ g)^*\nabla=f^*(g^*\nabla);
\end{equation}
moreover, if\/ $i:U\to M$ denotes the inclusion map of an open subset $U$ of $M$ then $i^*E$ can be naturally
identified with the bundle $E\vert_U$ and $i^*\nabla$ coincides with the induced connection $\nabla^U$.
\end{prop}
Identity \eqref{eq:pullpull} can be interpreted as a {\em chain rule\/} as follows; given a section $s:N\to E$
of $E$ along $f$ and $v\in TP$ then:
\[\big((f\circ g)^*\nabla\big)_v(s\circ g)\stackrel{\text{by \eqref{eq:pullpull}}}=\big(g^*(f^*\nabla)\big)_v(s\circ g)
\stackrel{\text{by \eqref{eq:defnablapull}}}=(f^*\nabla)_{\dd g(v)}s.\]

We have the following natural formula to compute the curvature and the torsion of a pull-back connection.
\begin{prop}\label{thm:RTpullback}
Given smooth manifolds $M$, $N$, a smooth vector bundle $E$ over $M$ endowed with a connection
$\nabla$ and a smooth map $f:N\to M$ then
the curvature tensor of $f^*\nabla$ is given by:
\[R^{f^*\nabla}_x(v,w)e=R^\nabla_{f(x)}\big(\dd f(x)v,\dd f(x)w\big)e,\]
for all $x\in N$, $v,w\in T_xN$, $e\in(f^*E)_x=E_{f(x)}$. Moreover, given a smooth vector bundle morphism
$\iota:TM\to E$, then $\iota\circ\dd f:TN\to E$ is identified with a vector bundle morphism $\tilde\iota:TN\to f^*E$
and the the following formula holds:
\begin{equation}\label{eq:tildeiota}
T^{\tilde\iota}_x(v,w)=T^\iota_{f(x)}\big(\dd f(x)v,\dd f(x)w\big),
\end{equation}
for all $x\in N$, $v,w\in T_xN$.
\end{prop}
Observe that if $E=TM$ and $\iota$ is the identity of $TM$ then formula \eqref{eq:defnablapull} means that:
\[(f^*\nabla)_X\big(\dd f(Y)\big)-(f^*\nabla)_Y\big(\dd f(X)\big)-\dd f\big([X,Y]\big)=T\big(\dd f(X),\dd f(Y)\big),\]
for all $X,Y\in\Gamma(TN)$.

Now we consider constructions acting on the fibers of the vector bundles. To this aim, we need some categorical language.
Given an integer $n\ge1$, we denote by $\Vect n$ the category whose objects are $n$-tuples $(V_i)_{i=1}^n$
of real finite-dimensional vector spaces and whose morphisms from $(V_i)_{i=1}^n$ to $(W_i)_{i=1}^n$ are $n$-tuples
$(T_i)_{i=1}^n$ of {\em vector space isomorphisms\/} $T_i:V_i\to W_i$. We set $\Vect 1=\Vect{}$.
A functor $\mathfrak F:\Vect n\to\Vect{}$ is called {\em smooth\/} if for any object $(V_i)_{i=1}^n$ of $\Vect n$
the map:
\begin{equation}\label{eq:frakFGL}
\mathfrak F:\GL(V_1)\times\cdots\times\GL(V_n)\longrightarrow\GL\big(\mathfrak F(V_1,\ldots,V_n)\big)
\end{equation}
is smooth. Observe that \eqref{eq:frakFGL} is a Lie group homomorphism; its differential at the identity is a
Lie algebra homomorphism that will be denoted by:
\[\overline{\mathfrak F}:\gl(V_1)\times\cdots\times\gl(V_n)\longrightarrow\gl\big(\mathfrak F(V_1,\ldots,V_n)\big).\]
Given vector bundles $E^1$, \dots, $E^n$ over a smooth manifold $M$ we obtain naturally a new vector bundle
$\mathfrak F(E^1,\ldots,E^n)$ over $M$ whose fiber at a point $x\in M$ is equal to $\mathfrak F(E^1_x,\ldots,E^n_x)$.
Given a smooth manifold $N$ and a smooth map $f:N\to M$, we may identify vector bundles
$f^*\big(\mathfrak F(E^1,\ldots,E^n)\big)$ and $\mathfrak F(f^*E^1,\ldots,f^*E^n)$.
Given vector bundle isomorphisms $L^i:E^i\to\widetilde E^i$, $i=1,\ldots,n$, then we obtain a vector bundle
isomorphism $L=\mathfrak F(T^1,\ldots,T^n)$ from $\mathfrak F(E^1,\ldots,E^n)$ to $\mathfrak F(\widetilde E^1,\ldots,\widetilde E^n)$
by setting:
\[L_x=\mathfrak F(L^1_x,\ldots,L^n_x),\]
for all $x\in M$.

We have the following functorial construction for connections.
\begin{prop}
Given an integer $n\ge1$ and a smooth functor $\mathfrak F:\Vect n\to\Vect{}$ then there exists a unique rule
that associates to each smooth manifold $M$, each $n$-tuple of vector bundles $(E^1,\ldots,E^n)$ over $M$
and each $n$-tuple of connections $(\nabla^1,\ldots,\nabla^n)$ on $(E^1,\ldots,E^n)$ respectively,
a connection $\nabla=\mathfrak F(\nabla^1,\ldots,\nabla^n)$ on $\mathfrak F(E^1,\ldots,E^n)$ satisfying the following
properties:
\begin{itemize}
\item[(a)] (naturality with pull-backs) given smooth manifolds $N$, $M$ and a smooth map $f:N\to M$ then $f^*\big(\mathfrak F(\nabla^1,\ldots,\nabla^n)\big)=
\mathfrak F(f^*\nabla^1,\ldots,f^*\nabla^n)$;
\item[(b)] (naturality with morphisms) given vector bundle isomorphisms $L^i:E^i\to\widetilde E^i$, $i=1,\ldots,n$, if $\nabla^i$
is a connection on $E^i$ which is $L^i$-related with a connection $\widetilde\nabla^i$ on $\widetilde E^i$ then
$\mathfrak F(\nabla^1,\ldots,\nabla^n)$ is $\mathfrak F(L^1,\ldots,L^n)$-related with
$\mathfrak F(\widetilde\nabla^1,\ldots,\widetilde\nabla^n)$;
\item[(c)] given connections $\nabla^i$ and $\widetilde\nabla^i$ on $E^i$ with
$\nabla^i-\widetilde\nabla^i=\mathfrak t^i$, $i=1,\ldots,n$, then:
\[\mathfrak F(\nabla^1,\ldots,\nabla^n)_Xs-\mathfrak F(\widetilde\nabla^1,\ldots,\widetilde\nabla^n)_Xs=
\overline{\mathfrak F}\big(\mathfrak t^1(X),\ldots,\mathfrak t^n(X)\big)s,\]
for all $s\in\Gamma\big(\mathfrak F(E^1,\ldots,E^n)\big)$;
\item[(d)] (trivial bundle property) If\/ $\nabla^i$ is the standard connection of the trivial
bundle $M\times E^i_0$ then $\mathfrak F(\nabla^1,\ldots,\nabla^n)$ is the standard connection
of the trivial bundle $M\times\mathfrak F(E^1_0,\ldots,E^n_0)$.
\end{itemize}
\end{prop}

Let $\mathfrak F=(\mathfrak F^1,\ldots,\mathfrak F^m)$ be an $m$-tuple of functors $\mathfrak F^i:\Vect n\to\Vect{}$
and let $\mathfrak G:\Vect m\to\Vect{}$ be a functor; we denote by $\mathfrak G\circ\mathfrak F:\Vect n\to\Vect{}$
the smooth functor defined\footnote{%
We will usually only describe functors on objects; the action of the functor on morphisms should be clear.} by:
\[(\mathfrak G\circ\mathfrak F)(V_1,\ldots,V_n)=\mathfrak G\big(\mathfrak F^1(V_1,\ldots,V_n),\ldots
\mathfrak F^m(V_1,\ldots,V_n)\big),\]
for all objects $V_1$, \dots, $V_n$ of $\Vect{}$.

\begin{prop}
Let $\mathfrak F=(\mathfrak F^1,\ldots,\mathfrak F^m)$ be an $m$-tuple of smooth functors $\mathfrak F^i:\Vect n\to\Vect{}$
and let $\mathfrak G:\Vect m\to\Vect{}$ be a smooth functor.
Given vector bundles $E^1$, \dots, $E^n$ over a smooth manifold $M$ endowed respectively
with connections $\nabla^1$, \dots, $\nabla^n$ then:
\[(\mathfrak G\circ\mathfrak F)(\nabla^1,\ldots,\nabla^n)=\mathfrak G\big(\mathfrak F^1(\nabla^1,\ldots,\nabla^n),\ldots,
\mathfrak F^m(\nabla^1,\ldots,\nabla^n)\big).\]
Moreover, if $\mathfrak I:\Vect{}\to\Vect{}$ denotes the identity functor of $\Vect{}$ then, given a connection
$\nabla$ on a vector bundle $E$, we have:
\[\mathfrak I(\nabla)=\nabla.\]
\end{prop}

\begin{prop}
Given a smooth functor $\mathfrak F:\Vect n\to\Vect{}$ and smooth vector bundles $\pi^i:E^i\to M$
endowed with connections $\nabla^i$, $i=1,\ldots,n$, then the curvature tensor of the connection
$\mathfrak F(\nabla^1,\ldots,\nabla^n)$ is given by:
\[R_x(v,w)=\overline{\mathfrak F}\big(R^1_x(v,w),\ldots,R^n_x(v,w)\big),\]
for all $x\in M$, $v,w\in T_xM$, where $R^i$ denotes the curvature tensor of $\nabla^i$, $i=1,\ldots,n$.
\end{prop}

\begin{defin}
Given a positive integer $n$ and smooth functors $\mathfrak F:\Vect n\to\Vect{}$ and $\mathfrak F':\Vect n\to\Vect{}$
then a {\em smooth natural transformation\/} $\rho$ from $\mathfrak F$ to $\mathfrak F'$ is a rule
that associates to each object $(V_i)_{i=1}^n$ of $\Vect n$ an open subset $A_{(V_1,\ldots,V_n)}$
of $\mathfrak F(V_1,\ldots,V_n)$ and a smooth map:
\[\rho_{V_1,\ldots,V_n}:A_{V_1,\ldots,V_n}\longrightarrow\mathfrak F'(V_1,\ldots,V_n)\]
such that, given objects $(V_i)_{i=1}^n$, $(W_i)_{i=1}^n$ of $\Vect n$ and a morphism
$(T_i)_{i=1}^n$ from $(V_i)_{i=1}^n$ to $(W_i)_{i=1}^n$ then $\mathfrak F(T_1,\ldots,T_n)$ carries
$A_{V_1,\ldots,V_n}$ to $A_{W_1,\ldots,W_n}$ and the diagram:
\[\xymatrix@C+25pt{A_{V_1,\ldots,V_n}\ar[r]^-{\rho_{V_1,\ldots,V_n}}\ar[d]_{\mathfrak F(T_1,\ldots,T_n)}&
\mathfrak F'(V_1,\ldots,V_n)\ar[d]^{\mathfrak F'(T_1,\ldots,T_n)}\\
A_{W_1,\ldots,W_n}\ar[r]_-{\rho_{W_1,\ldots,W_n}}&\mathfrak F'(W_1,\ldots,W_n)}\]
commutes.
\end{defin}
Given a smooth natural transformation $\rho$ from $\mathfrak F$ to $\mathfrak F'$ and given vector bundles
$\pi^i:E^i\to M$, $i=1,\ldots,n$, we obtain a fiber bundle morphism $\rho_{E^1,\ldots,E^n}:A\to\mathfrak F'(E^1,\ldots,E^n)$
defined on an open subset $A$ of $\mathfrak F(E^1,\ldots,E^n)$ by setting:
\[A_x=A_{E^1_x,\ldots,E^n_x},\quad(\rho_{E^1,\ldots,E^n})_x=\rho_{E^1_x,\ldots,E^n_x},\]
for all $x\in M$.

\begin{prop}\label{thm:connectionNT}
Given a positive integer $n$, smooth functors $\mathfrak F:\Vect n\to\Vect{}$, $\mathfrak F':\Vect n\to\Vect{}$,
a smooth natural transformation $\rho$ from $\mathfrak F$ to $\mathfrak F'$, vector bundles $\pi^i:E^i\to M$
endowed with connections $\nabla^i$, $i=1,\ldots,n$, then:
\[\mathfrak F'(\nabla^1,\ldots,\nabla^n)_v(\rho_{E^1,\ldots,E^n}\circ s)=\dd\rho_{E^1_x,\ldots,E^n_x}\big(s(x)\big)
\big(\mathfrak F(\nabla^1,\ldots,\nabla^n)_vs\big),\]
for all $x\in M$, $v\in T_xM$ and every smooth section $s$ of $\mathfrak F(E^1,\ldots,E^n)$
with range contained in the domain of $\rho_{E^1,\ldots,E^n}$.
\end{prop}

\begin{example}
Let $n$ be a positive integer and consider the smooth functor $\mathfrak S:\Vect n\to\Vect{}$ defined by:
\[\mathfrak S(V_1,\ldots,V_n)=V_1\oplus\cdots\oplus V_n.\]
Given vector bundles $E^1$, \dots, $E^n$ over a smooth manifold $M$ then $\mathfrak S(E^1,\ldots,E^n)$
is the Whitney sum of $E^1$, \dots, $E^n$. Let $\nabla^i$ be a connection on $E^i$, $i=1,\ldots,n$.
For each $i=1,\ldots,n$, consider the smooth functor $\mathfrak P^i:\Vect n\to\Vect{}$ defined by:
\[\mathfrak P^i(V_1,\ldots,V_n)=V_i.\]
We have a smooth natural transformation $\rho^i$ from $\mathfrak S$ to $\mathfrak P^i$ given by:
\[\rho^i:V_1\oplus\cdots\oplus V_n\ni(v_1,\ldots,v_n)\longmapsto v_i\in V_i.\]
Set $\nabla=\mathfrak S(\nabla^1,\ldots,\nabla^n)$.
Proposition~\ref{thm:connectionNT} implies that:
\[\nabla_v(s_1,\ldots,s_n)=(\nabla^1_vs_1,\ldots,\nabla^n_vs_n),\]
for all $s_1\in\Gamma(E^1)$, \dots, $s_n\in\Gamma(E^n)$, $v\in TM$.
\end{example}

\begin{example}
Consider the smooth functors $\mathfrak F:\Vect2\to\Vect{}$ and $\mathfrak G:\Vect2\to\Vect{}$ defined as follows;
let $V_1$, $V_2$, $W_1$, $W_2$ be objects of $\Vect{}$ and let $T_1:V_1\to W_1$, $T_2:V_2\to W_2$ be isomorphisms.
We set:
\begin{gather*}
\mathfrak F(V_1,V_2)=\Lin(V_1,V_2),\quad\mathfrak F(T_1,T_2)L=T_2\circ L\circ T_1^{-1},\\
\mathfrak G(V_1,V_2)=\Lin(V_2^*,V_1^*)\quad\mathfrak G(T_1,T_2)R=(T_1^*)^{-1}\circ R^*\circ T_2^*,
\end{gather*}
for $L\in\Lin(V_1,V_2)$, $R\in\Lin(V_2^*,V_1^*)$. We have a natural transformation $\rho$ from $\mathfrak F$
to $\mathfrak G$ defined by:
\[\rho:\Lin(V_1,V_2)\ni t\longmapsto t^*\in\Lin(V_2^*,V_1^*).\]
Let $E^1$, $E^2$ be vector bundles over a smooth manifold $M$ endowed with connections $\nabla^1$ and $\nabla^2$
respectively. We denote by $\nabla$ both the connections $\mathfrak F(\nabla^1,\nabla^2)$ and $\mathfrak G(\nabla^1,\nabla^2)$
on the bundles $\Lin(E^1,E^2)$ and $\Lin((E^2)^*,(E^1)^*)$ respectively. Proposition~\ref{thm:connectionNT}
tells us that, given a smooth section $L$ of $\Lin(E^1,E^2)$ then:
\[\nabla_vL^*=(\nabla_vL)^*,\]
for all $v\in TM$.
\end{example}

\begin{example}
Consider the smooth functor $\mathfrak F:\Vect{}\to\Vect{}$ defined as follows;
let $V$, $W$ be objects of $\Vect{}$ and let $T:V\to W$ be an isomorphism. We set:
\begin{gather*}
\mathfrak F(V)=\Bilin(V,V;V)\oplus V\oplus V,\\
\mathfrak F(T)(B,v_1,v_2)=\big(T\circ B(T^{-1}\cdot,T^{-1}\cdot),T(v_1),T(v_2)\big),
\end{gather*}
for every bilinear map $B:V\times V\to V$ and all $v_1,v_2\in V$. We have a smooth natural transformation
$\rho$ from $\mathfrak F$ to the identity functor $\mathfrak I$ of $\Vect{}$ defined by:
\[\rho:\Bilin(V,V;V)\oplus V\oplus V\ni(B,v_1,v_2)\longmapsto B(v_1,v_2)\in V.\]
Let $E$ be a vector bundle over a smooth manifold $M$ endowed with a connection $\nabla$.
We will also denote by $\nabla$ the connection $\mathfrak F(\nabla)$ on $\Bilin(E,E;E)$.
Proposition~\ref{thm:connectionNT} tells us that, given a smooth section $B$ of $\Bilin(E,E;E)$
and smooth sections $s_1$, $s_2$ of $E$ then:
\[\nabla_v\big(B(s_1,s_2)\big)=(\nabla_vB)(s_1,s_2)+B(\nabla_vs_1,s_2)+B(s_1,\nabla_vs_2),\]
for all $v\in TM$.
\end{example}

\begin{example}
Consider the smooth functor $\mathfrak F:\Vect{}\to\Vect{}$ defined as follows;
let $V$, $W$ be objects of $\Vect{}$ and let $T:V\to W$ be an isomorphism. We set:
\begin{gather*}
\mathfrak F(V)=\Lin(V),\\
\mathfrak F(T)L=T\circ L\circ T^{-1},
\end{gather*}
for all $L\in\Lin(V)$. We have a smooth natural transformation $\rho$ from $\mathfrak F$ to itself defined by:
\[\rho:\Lin(V)\supset\GL(V)\ni L\longmapsto L^{-1}\in\Lin(V).\]
Let $E$ be a vector bundle over a smooth manifold $M$ endowed with a connection $\nabla$.
We will also denote by $\nabla$ the connection $\mathfrak F(\nabla)$ on $\Lin(E)$.
Let $L$ be a smooth section of $\Lin(E)$ such that $L_x$ is an isomorphism of $E_x$, for all $x\in M$.
Proposition~\ref{thm:connectionNT} tells us that:
\[\nabla_v(L^{-1})=-L^{-1}(\nabla_vL)L^{-1},\]
for all $v\in TM$.
\end{example}

\end{section}


\begin{thebibliography}{99}

\bibitem{KN} S. Kobayashi and K. Nomizu, {\em Foundations of differential geometry. Vol I},
Interscience Publishers, New York-London, 1963.

\bibitem{Lang} S.\ Lang, \emph{Fundamentals of differential geometry},
Graduate Texts in Mathematics. 191. New York, NY: Springer, 2001.

\bibitem{Cordoba} P. Piccione, D. V. Tausk, \emph{Connections compatible with tensors.
A characterization of left-invariant Levi--Civita connections in Lie groups}, arXiv
\texttt{math.DG/0509656}.

\bibitem{Wolf} J. A. Wolf, {\em Spaces of constant curvature}. New York, Mc-Graw-Hill, 1967.
\end{thebibliography}
\end{document}